\documentclass[11pt]{article}

\usepackage[ansinew]{inputenc}
\usepackage{amsfonts}
\usepackage{amssymb,amsmath}
\usepackage{latexsym}
\usepackage{graphics}
\usepackage{epic}
\usepackage{epsfig}
\usepackage{psfrag}
\usepackage{subfigure}
\usepackage{bbm}

\numberwithin{equation}{section}

\newcommand     {\ZZ}{\mathbb{Z}}
\newcommand     {\NN}{\mathbb{N}}
\newcommand     {\RR}{\mathbb{R}}
\newcommand     {\PP}{\mathbb{P}}

\newcommand     {\EE}{\mathbb{E}}

\newtheorem     {thm}{Theorem}[section]

\newtheorem     {lem}[thm]{Lemma}
\newtheorem     {prop}[thm]{Proposition}
\newtheorem     {cor}[thm]{Corollary}

\begin{document}

\title{Splitting trees with neutral Poissonian mutations II: Largest and oldest
  families.}

\author{\textsc{Nicolas Champagnat$^{1}$, Amaury Lambert$^{2}$}}

\footnotetext[1]{TOSCA project-team, INRIA Nancy -- Grand Est, IECN -- UMR 7502,
  Nancy-Universit\'e, Campus scientifique, B.P.\ 70239, 54506 Vand\oe uvre-l\`es-Nancy Cedex,
  France, E-mail: \texttt{Nicolas.Champagnat@inria.fr}}

\footnotetext[2]{Laboratoire de Probabilités et Modèles Aléatoires,
UMR 7599 CNRS and UPMC Univ Paris 06,
Case courrier 188, 
4 Place Jussieu,
F-75252 Paris Cedex 05, France, Email: \texttt{amaury.lambert@upmc.fr}
}
\maketitle

\begin{abstract}
We consider a supercritical branching population, where individuals have i.i.d.\ lifetime durations (which are not necessarily exponentially distributed) and give birth (singly) at constant rate.
We assume that individuals independently experience neutral mutations, at constant rate $\theta$ during their lifetimes, under the infinite-alleles assumption: each mutation instantaneously confers a brand new type, called allele or haplotype, to its carrier. The type carried by a mother at the time when she gives birth is transmitted to the newborn.

We are interested in the sizes and ages at time $t$ of the clonal families carrying the most abundant alleles or the oldest ones, as $t\to\infty$, on the survival event.  Intuitively, the results must depend on how the mutation rate $\theta$ and the Malthusian parameter $\alpha>0$ compare. Hereafter, $N\equiv N_t$ is the population size at time $t$, constants $a,c$ are scaling constants, whereas $k,k'$ are explicit positive constants which depend on the parameters of the model.

When $\alpha >\theta$, 
the most abundant families are also the oldest ones, they have size $cN^{1-\theta/\alpha}$
 and age $t-a$. 

When $\alpha<\theta$, the oldest families have age $(\alpha /\theta)t + a$ and tight sizes; the most abundant families have sizes 
$k\log(N) -k' \log \log (N)+c$
and all have age $(\theta- \alpha)^{-1} \log(t)$. 

When $\alpha=\theta$, the oldest families have age $kt -k'\log(t)+ a$ and tight sizes; the most abundant families have sizes $(k\log(N) -k' \log \log (N)+c)^2$
and all have age $t/2$.

Those informal results can be stated rigorously in expectation. Relying heavily on the theory of coalescent point processes \cite{L10, P}, we are also able, when $\alpha \le \theta$, to show convergence in distribution of the joint, properly scaled ages and sizes of the most abundant/oldest families and to specify the limits as some explicit Cox processes.

This is in deep contrast with the largest/oldest families in the standard coalescent with Poissonian mutations, which converge to some point processes after being rescaled by $N$ \cite{DT, D, Ewens}.
\end{abstract}          
\bigskip

\bigskip


\section{Introduction and motivation}
\label{sec:intro}

We consider a general branching population, where individuals reproduce independently of each
other, have i.i.d.\ lifetime durations with arbitrary distribution, and give birth at constant rate during their lifetime. We also assume that each birth gives rise to a single newborn. The genealogical tree associated with this construction is known as a splitting tree \cite{Geiger, GK, L10}.
The process $(N_t;t\ge 0)$ counting the population size is a non-Markovian birth--death process belonging to the class of general branching processes, or \emph{Crump--Mode--Jagers} (CMJ) processes. Since births arrive singly and at constant rate, these processes are sometimes called \emph{homogeneous, binary CMJ processes}. 

Also, individuals are given a type, called allele or haplotype. They inherit their type at birth from their mother, and (their germ line) can change type throughout their lifetime, at the points of independent Poisson point processes with rate $\theta$, conditional on lifetimes (neutral mutations). The type conferred by a mutation is each time an entirely new type, an assumption known as the \emph{infinitely-many alleles model}.

We are interested in the so-called \emph{allelic partition} (partition into types) of the population alive at time
$t$. In~\cite{CL1, L09}, we obtained explicit formulae for the expected \emph{frequency spectrum} of the allelic
partition. The frequency spectrum is a convenient way of describing this partition without labelling types. It
is defined as the sequence of numbers $(A_\theta(k,t),k\geq 0)$, where $A_\theta(k,t)$ is the number of types
carried by $k$ individuals at time $t$.  For example in \cite{CL1}, we have derived explicit formulae for the expectation of $A_\theta(k,t)$ conditional on population size at $t$. From these formulae, using the theory of branching processes counted by random characteristics, we have specified the a.s. limit, as $t\to\infty$, of the fraction of types carried by a fixed number $k$ of individuals.  

If we call \emph{clonal families}, or simply \emph{families}, the equivalence classes associated to identity by type (i.e., the components of the allelic partition), it is usual to call \emph{small families} the families of sizes $k=1,2,3,$... Here, \emph{large families} will refer to families with most frequent (i.e., abundant) types having alive representatives, and \emph{old families} to families having oldest types with alive representatives, where the age of a type is the time elapsed since its original mutation.
In the present work, we are interested in the asymptotic behavior, as $t\to\infty$, of the
sizes and ages of large and of old families. 

The most celebrated mathematical result regarding allelic partitions is Ewens' sampling formula, which provides the distribution of the frequency spectrum for the Kingman coalescent tree with neutral Poissonian mutations \cite{Ewens}. It has notably been shown \cite{DT, D} that under this model, the largest (resp. oldest) families converge, after being rescaled by the population size $N$, to the Poisson--Dirichlet (resp. GEM) distribution. We will see here, that, for example, the largest families are never of the order of $N$, but depending on how the Malthusian parameter $\alpha$ scales with the mutation rate $\theta$, of the order of $N^{1-\theta/\alpha}$ (case $\theta<\alpha$), order of $(\log N)^2$ (case $\theta=\alpha$), or order of $\log N$ (case $\theta>\alpha$). The first case ($\theta<\alpha$) shows more similarities with the frequency spectrum of the Beta-coalescent~\cite{BBS}.

We refer the reader to \cite{CL1} for more references on the topic of allelic partitions, especially regarding branching processes. For example, similar questions were  studied for general CMJ processes, when mutations occur
\emph{at birth}, in the monography due to Z.\ Ta\"ib~\cite{Taib}. These results
rely heavily on the theory of branching processes counted by random characteristics, due to
P.\ Jagers and O.\ Nerman \cite{J, JNa, JNb, N}, and more specifically on \emph{time dependent
  random characteristics} as developed in~\cite{JNb}. Z.\ Ta\"ib obtains results of
convergence in distribution of the (correctly rescaled) point process of ages, similar to the
results we obtain in Sections~\ref{sec:subcrit} and~\ref{sec:crit}. However, the techniques
of~\cite{Taib} do not apply to the case where mutations occur during individuals' lifetimes,
since the genealogical tree of types is not the one of a CMJ process in this case.

One of the initial motivations of \cite{CL1} and of the present work was the following model inspired from the works of P.C.\ Sabeti and her coauthors (see e.g.~\cite{Sabeti}). Inside a large population, consider a subpopulation consisting of individuals carrying a specific selective gene called `core haplotype' and thus experiencing demographic growth. The haplotypic structure of the subpopulation restricted to a portion of length $L$ around the core haplotype on the chromosome carrying it, is assumed to be altered by recombination. As long as the total population is sufficiently large w.r.t. the growing subpopulation, each time a sequence belonging to (an individual in) this subpopulation recombines with another sequence, with high probability this sequence will be a new sequence belonging to the rest of the population. Therefore, the new sequence obtained after recombination can be treated as a mutant under the infinitely-many alleles model. In this setting, mutation rate is an increasing function of $L$. In \cite{Sabeti}, a tree representation of the allelic partition as a function of $L$ is given for each core haplotype in a given set of genes suspected to have been selected in humans. The tree obtained this way is called ``recombination tree''. An interesting question is to develop statistical methods allowing to detect positive selection from the knowledge of this tree. Here, we assume that our (sub)population grows at a Malthusian rate $\alpha$ (supercritical CMJ process). We are able to give the asymptotic distribution of the rightmost part of the frequency spectrum for a given mutation rate (this corresponds to fixing $L$ in the recombination setting). Since $\theta$ can be seen as a death rate when restricting the count to individuals carrying the same allele, the phase transition at $\theta=\alpha$ is intuitive. In the recombination tree, this phase transition should translate into a transition, at a certain locus length $L_0$, from a small number of thick branches ($L<L_0$) to a large number of thin branches ($L>L_0$). We plan to extend this study to a full description of the structure of the recombination tree.

In the next section, we define rigorously the model and recall some chosen results from \cite{CL1}. Section \ref{sec:expected} is concerned with the asymptotic behavior, as $t\to\infty$, of the expected sizes and ages of the most abundant/oldest families. Sections \ref{sec:subcrit} and \ref{sec:crit} deal with the joint convergence in distribution of these sizes and ages in the respective cases when clonal families are subcritical or critical. A final appendix is devoted to some technical lemmas for the control of moments of order 2 of largest sizes and ages.

\section{Model, preliminary results and statement of the main results}
\label{sec:models}

\subsection{Model}

In this work, we consider genealogical trees satisfying the branching property and called \emph{splitting trees} \cite{Geiger, GK}. Splitting trees are those random trees where individuals' lifetime durations are i.i.d.\ with an arbitrary distribution, but where birth events occur at Poisson times during each individual's lifetime. We call $b$ this constant birth rate and we denote by $V$ a r.v.\ distributed as the lifetime duration. Then set $\Lambda(dr):=b\PP(V\in dr)$ a finite measure on $(0,\infty]$ with total mass $b$ called the \emph{lifespan measure}. We will always assume that a splitting tree is started with one unique progenitor born at time 0.

The process $(N_t;t\ge0)$ counting the number of alive individuals at time $t$ is a homogeneous, binary \emph{Crump--Mode--Jagers process}, which is not Markovian unless $\Lambda$ has an exponential density or is a Dirac mass at $\{+\infty\}$.

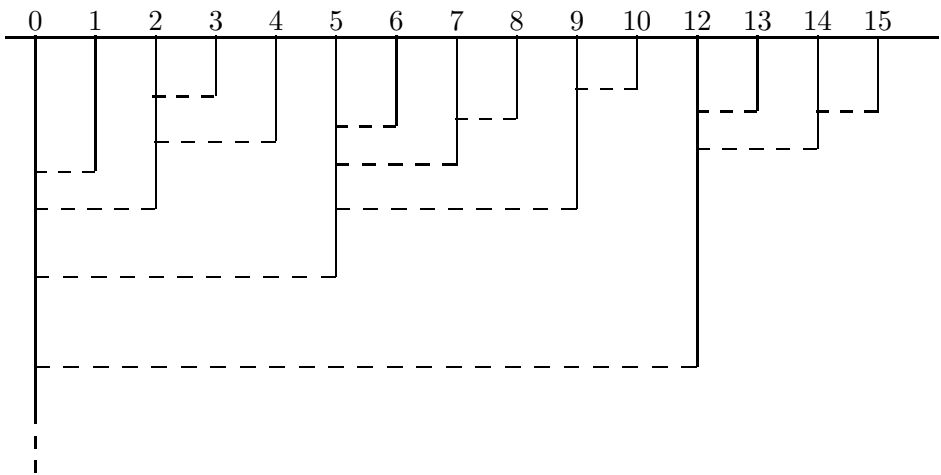
\begin{figure}[ht]

\unitlength 2mm 
\linethickness{0.4pt}

\begin{picture}(66,33)(-5,10)
\put(4,39.875){\line(1,0){62}}
\put(10,40){\line(0,-1){9}}
\put(14,40){\line(0,-1){11.5}}
\put(18,40){\line(0,-1){4}}
\put(22,40){\line(0,-1){7}}
\put(26,40){\line(0,-1){16}}
\put(30,40){\line(0,-1){6}}
\put(34,39.875){\line(0,-1){8.5}}
\put(38,40){\line(0,-1){5.5}}
\put(42,40){\line(0,-1){11.5}}
\put(46,40){\line(0,-1){3.625}}
\put(50,40){\line(0,-1){22}}
\put(54,40){\line(0,-1){5}}
\put(58,40){\line(0,-1){7.5}}
\put(62,39.875){\line(0,-1){5}}
\put(6,40){\line(0,-1){25}}
\put(6,15){\line(0,-1){.8}}
\put(6,13.33){\line(0,-1){.8}}
\put(6,11.73){\line(0,-1){.8}}
\put(9.93,30.93){\line(-1,0){.8}}
\put(8.33,30.93){\line(-1,0){.8}}
\put(6.73,30.93){\line(-1,0){.8}}
\put(13.93,28.43){\line(-1,0){.8889}}
\put(12.152,28.43){\line(-1,0){.8889}}
\put(10.374,28.43){\line(-1,0){.8889}}
\put(8.596,28.43){\line(-1,0){.8889}}
\put(6.819,28.43){\line(-1,0){.8889}}
\put(17.805,35.93){\line(-1,0){.8}}
\put(16.205,35.93){\line(-1,0){.8}}
\put(14.605,35.93){\line(-1,0){.8}}
\put(21.93,32.93){\line(-1,0){.8889}}
\put(20.152,32.93){\line(-1,0){.8889}}
\put(18.374,32.93){\line(-1,0){.8889}}
\put(16.596,32.93){\line(-1,0){.8889}}
\put(14.819,32.93){\line(-1,0){.8889}}
\put(25.93,23.93){\line(-1,0){.9524}}
\put(24.025,23.93){\line(-1,0){.9524}}
\put(22.12,23.93){\line(-1,0){.9524}}
\put(20.215,23.93){\line(-1,0){.9524}}
\put(18.311,23.93){\line(-1,0){.9524}}
\put(16.406,23.93){\line(-1,0){.9524}}
\put(14.501,23.93){\line(-1,0){.9524}}
\put(12.596,23.93){\line(-1,0){.9524}}
\put(10.692,23.93){\line(-1,0){.9524}}
\put(8.787,23.93){\line(-1,0){.9524}}
\put(6.882,23.93){\line(-1,0){.9524}}
\put(29.93,33.93){\line(-1,0){.8}}
\put(28.33,33.93){\line(-1,0){.8}}
\put(26.73,33.93){\line(-1,0){.8}}
\put(33.93,31.43){\line(-1,0){.8889}}
\put(32.152,31.43){\line(-1,0){.8889}}
\put(30.374,31.43){\line(-1,0){.8889}}
\put(28.596,31.43){\line(-1,0){.8889}}
\put(26.819,31.43){\line(-1,0){.8889}}
\put(37.93,34.43){\line(-1,0){.8}}
\put(36.33,34.43){\line(-1,0){.8}}
\put(34.73,34.43){\line(-1,0){.8}}
\put(41.93,28.43){\line(-1,0){.9412}}
\put(40.047,28.43){\line(-1,0){.9412}}
\put(38.165,28.43){\line(-1,0){.9412}}
\put(36.283,28.43){\line(-1,0){.9412}}
\put(34.4,28.43){\line(-1,0){.9412}}
\put(32.518,28.43){\line(-1,0){.9412}}
\put(30.636,28.43){\line(-1,0){.9412}}
\put(28.753,28.43){\line(-1,0){.9412}}
\put(26.871,28.43){\line(-1,0){.9412}}
\put(45.93,36.43){\line(-1,0){.8}}
\put(44.33,36.43){\line(-1,0){.8}}
\put(42.73,36.43){\line(-1,0){.8}}
\put(49.93,17.93){\line(-1,0){.9778}}
\put(47.974,17.93){\line(-1,0){.9778}}
\put(46.019,17.93){\line(-1,0){.9778}}
\put(44.063,17.93){\line(-1,0){.9778}}
\put(42.107,17.93){\line(-1,0){.9778}}
\put(40.152,17.93){\line(-1,0){.9778}}
\put(38.196,17.93){\line(-1,0){.9778}}
\put(36.241,17.93){\line(-1,0){.9778}}
\put(34.285,17.93){\line(-1,0){.9778}}
\put(32.33,17.93){\line(-1,0){.9778}}
\put(30.374,17.93){\line(-1,0){.9778}}
\put(28.419,17.93){\line(-1,0){.9778}}
\put(26.463,17.93){\line(-1,0){.9778}}
\put(24.507,17.93){\line(-1,0){.9778}}
\put(22.552,17.93){\line(-1,0){.9778}}
\put(20.596,17.93){\line(-1,0){.9778}}
\put(18.641,17.93){\line(-1,0){.9778}}
\put(16.685,17.93){\line(-1,0){.9778}}
\put(14.73,17.93){\line(-1,0){.9778}}
\put(12.774,17.93){\line(-1,0){.9778}}
\put(10.819,17.93){\line(-1,0){.9778}}
\put(8.863,17.93){\line(-1,0){.9778}}
\put(6.907,17.93){\line(-1,0){.9778}}
\put(53.93,34.93){\line(-1,0){.8}}
\put(52.33,34.93){\line(-1,0){.8}}
\put(50.73,34.93){\line(-1,0){.8}}
\put(57.93,32.43){\line(-1,0){.8889}}
\put(56.152,32.43){\line(-1,0){.8889}}
\put(54.374,32.43){\line(-1,0){.8889}}
\put(52.596,32.43){\line(-1,0){.8889}}
\put(50.819,32.43){\line(-1,0){.8889}}
\put(61.93,34.93){\line(-1,0){.8}}
\put(60.33,34.93){\line(-1,0){.8}}
\put(58.73,34.93){\line(-1,0){.8}}
\put(6,41){\makebox(0,0)[cc]{$0$}}
\put(10,41){\makebox(0,0)[cc]{$1$}}
\put(14,41){\makebox(0,0)[cc]{$2$}}
\put(18,41){\makebox(0,0)[cc]{$3$}}
\put(22,41){\makebox(0,0)[cc]{$4$}}
\put(26,41){\makebox(0,0)[cc]{$5$}}
\put(30,41){\makebox(0,0)[cc]{$6$}}
\put(34,41){\makebox(0,0)[cc]{$7$}}
\put(38,41){\makebox(0,0)[cc]{$8$}}
\put(42,41){\makebox(0,0)[cc]{$9$}}
\put(46,41){\makebox(0,0)[cc]{$10$}}
\put(50,41){\makebox(0,0)[cc]{$12$}}
\put(54,41){\makebox(0,0)[cc]{$13$}}
\put(58,41){\makebox(0,0)[cc]{$14$}}
\put(62,41){\makebox(0,0)[cc]{$15$}}
\end{picture}

\caption{A coalescent point process.}
\label{fig : coalpointproc}
\end{figure}

In \cite{L10}, it is shown that the genealogy of a splitting tree conditioned to be extant at a fixed time  $t$ is given by a \emph{coalescent point process}, that is, a sequence of i.i.d.\ random variables $H_i$, $i\ge 1$, killed at its first value greater than $t$. In particular, conditional on $N_t\not=0$, $N_t$ follows a geometric ditribution with parameter $\PP(H<t)$. More specifically, for any $0\le i \le N_t-1$, the \emph{coalescence time}  between the $i$-th individual alive at time $t$ and the $j$-th individual alive at time $t$ (i.e., the time elapsed since the common lineage to both individuals split into two distinct lineages) is the maximum of $H_{i+1},\ldots, H_j$. The graphical representation on Figure \ref{fig : coalpointproc} is straightforward. The common law of these so-called \emph{branch lengths} is given by
\begin{equation}
  \label{eq:def-law H}
\PP(H>s)=\frac{1}{W(s)} , 
\end{equation}
where the nondecreasing function $W$ is such that $W(0)=1$ and is characterized by its Laplace transform. More
specifically, these branch lengths are the depths of the excursions of the jump contour process, say $Y^{(t)}$, of the splitting tree truncated below level $t$. They are i.i.d.\ because $Y^{(t)}$ is a Markov process. Indeed, it is shown in \cite{L10} that $Y^{(t)}$ has the law of a Lévy process, say $Y$, with no negative jumps, reflected below $t$ and killed upon hitting 0. The function $W$ is called the \emph{scale function} of $Y$, and is defined from the Laplace exponent $\psi$ of $Y$:
\begin{equation}
  \label{eq:def-psi}
  \psi(x)=x-\int_{(0,+\infty]}\left(1-e^{-rx}\right) \Lambda (dr)\qquad x\in \RR_+.  
\end{equation}
Let $\alpha$ denote the largest root of $\psi$. In the supercritical case (i.e.\ $\int_{(0,\infty]} r\Lambda (dr)>1$), and in this case only, $\alpha$ is positive and called the \emph{Malthusian parameter}, because the population size grows exponentially at rate $\alpha$ on the survival event. Then the function $W$ is characterized by
$$
\int_0^\infty e^{-xr} W(r)\, dr = \frac{1}{\psi(x)}\qquad x>\alpha.
$$
Actually, it is possible to show by path decompositions of the process $Y$ that 
\begin{equation}
  \label{eq:link-W-excursion}
  W(x)=\exp\left( b\int_0^x dt\,\PP(J>t) \right),  
\end{equation}
where $J$ is the maximum of the path of $Y$ killed upon hitting 0 and started from a random initial value, distributed as $V$. Note that since $Y$ is also the contour process of a splitting tree, $J$ has the law of the extinction time of the CMJ process $N$ started from one individual.


Throughout this work, we further assume that individuals independently experience mutations at Poisson times during
their lifetime, that each new mutation event confers a brand new type (called haplotype, or allele) to the
individual, and that a newborn holds the same type as her mother at birth time. The mutation rate is denoted by $\theta$. From now on, $\PP_t$ (resp. $\PP^\star$) will denote the conditional probability on survival up to time $t$ (resp. on the survival event).

\subsection{Expected frequency spectrum}
Recall from the introduction that $A_\theta(k,t)$ is the number of types carried by $k$ individuals at time $t$. Also denote by $Z_0(t)$ the number of individuals carrying the ancestral type at time $t$.

In~\cite[Cor.\:4.3]{CL1}, we obtained the following explicit formulae for the expected frequency spectrum in
the population at time $t$:
\begin{equation}
  \label{eq:exp-FS}
  \EE_t A_\theta (k,t) =W(t)\int_0^t dx\,\theta\,e^{-\theta x}\, \frac{1}{W_\theta(x)^2}\left(1-\frac{1}{W_\theta(x)}\right)^{k-1},  
\end{equation}
and
\begin{equation}
  \label{eq:ancestral-FS}
  \PP_t\left(Z_0(t)=k\right)= W(t)\,\frac{e^{-\theta t}}{W_\theta(t)^2}\left(1-\frac{1}{W_\theta(t)}\right)^{k-1},  
\end{equation}
 and
\begin{equation}
  \label{eq:W_theta}
  W_\theta(x):=e^{-\theta x}W(x)+\theta\int_0^x W(y)e^{-\theta y}dy,
\end{equation}
i.e.\ $W'_\theta(x)=e^{-\theta x}W'(x)\geq 0$ and $W_\theta(0)=W(0)=1$. Note that $W_\theta$
is the scale function associated to the clonal splitting tree~\cite[Thm.\:3.1]{CL1}, i.e.
\begin{equation*}
  \int_0^\infty e^{-xr}W_\theta(r)dr=\frac{1}{\psi_\theta(x)},  \qquad x>\alpha-\theta,
\end{equation*}
where
\begin{equation}
  \label{eq:def-psi-theta}
  \psi_\theta(x):=x-\int_{(0,+\infty]}(1-e^{-rx})b\PP(V_\theta\in dr)=\frac{x\psi(x+\theta)}{x+\theta},
\end{equation}
where $V_\theta$ denotes the minimum of $V$ and of an independent r.v.\ with parameter
$\theta$.

In addition, we were able to obtain the explicit expected age density of the frequency
spectrum~\cite[Eq.\:(4.5)]{CL1}: defining $A_\theta(k,t,y;dy)$ as the number of haplotypes of age in the
interval $(y,y+dy)$ represented by exactly $k$ alive individuals at time $t$, we have
\begin{equation}
  \label{eq:exp-FS-with-age}
  \EE_t A_\theta (k,t,y;dy) =\theta\, dy\, W(t)\frac{e^{-\theta y}}{W_\theta(y)^2}\left(1-\frac{1}{W_\theta(y)}\right)^{k-1} .
\end{equation}

In~\cite{CL1}, denoting by $A_\theta(t)=\sum_{k\geq 1}A_\theta(k,t)$ the total number of haplotypes at time $t$, we deduced from these expressions the a.s.\ large time convergence of the fraction $A_\theta(k,t)/A_\theta(t)$. Recall that families with given size $k=1,2,3,...$ are referred to as  \emph{small families}. \emph{Large families} are those who have the largest sizes and \emph{old families} are those whose original mutation is among the oldest. We are interested in the sizes and ages of large and of old families. For example, the size of the largest family is the largest $k$ such that $A_\theta(k,t)\ge 1$.

\subsection{Statement of results}

Recall that we always assume that the Malthusian parameter $\alpha$ is positive. The
asymptotic size of the most frequent and oldest haplotypes strongly depends on the way
$\alpha$ and the mutation rate $\theta$ compare. Since $\theta$ is an additional death rate
for clonal families, the case $\alpha>\theta$ corresponds to supercritical clonal families,
the case $\theta=\alpha$ corresponds to critical clonal families, and the case $\theta>\alpha$
corresponds to subcritical clonal families.

In the whole paper, we are going to use the following notation: for all $x,s>0$, define
\begin{itemize}
\item $L_t(x)$ the number of haplotypes carried by more than (or exactly) $x$ individuals alive at time $t$
  ($L$ for \emph{large});
\item $O_t(s)$ the number of haplotypes with alive representatives at time $t$ older than $s$, i.e.\ whose original mutation has age greater than $s$ ($O$ for \emph{old});
\item $M_t(x,s)$ the number of haplotypes carried by more than $x$ individuals alive at time
  $t$ and whose original mutation has age greater than $s$. By convention, we set
  $M_t(x,s)=M_t(x,0)=L_t(x)$ when $s<0$. For $0\leq s_1\leq s_2$, we also define
  $$
  M_t(x,s_1,s_2)=M_t(x,s_1)-M_t(x,s_2),
  $$
  the number of haplotypes carried by more than $x$ individuals alive at time $t$, whose original mutation has age in  $(s_1, s_2]$.
\end{itemize}

Our convergence results are of two kinds: convergence in expectation of $L_t(x_t)$ and
$O_t(s_t)$ for conveniently chosen $x_t$ and $s_t$, which are directly obtained
from~(\ref{eq:exp-FS}),~(\ref{eq:ancestral-FS}) and~(\ref{eq:exp-FS-with-age}) (see
Section~\ref{sec:expected}), and convergence in distribution of the point process of
(correctly rescaled) largest families or oldest families, which require to combine the
previous equations with the theory of coalescent point processes~\cite{L10, P} (see
Sections~\ref{sec:subcrit} and~\ref{sec:crit}). We obtain different results depending on
whether the clonal families are supercritical, subcritical or critical.


\subsubsection{Supercritical clonal families}
\label{sec:results-supercrit}

Our only result in the supercritical case ($\alpha>\theta$) is the following
(Proposition~\ref{prop:frequent-a>t}): for all $0\leq a<b\leq+\infty$ and $c\geq 0$,
\begin{equation}
  \label{eq:result-supercrit}
  \lim_{t\rightarrow +\infty}\EE_t M_t(ce^{(\alpha-\theta)t},t-a_1,t-a_0)
  =\frac{\alpha-\theta}{\alpha}\int_{a_0}^{a_1}\exp\left(\alpha
    y-c\psi'_\theta(\alpha-\theta)e^{(\alpha-\theta)y}\right)\,(\theta dy+\delta_0(dy)),
\end{equation}
where $\delta_0$ denotes the Dirac measure at 0. Note that~(\ref{eq:def-psi-theta}) yields
\begin{equation*}
  \psi'_\theta(\alpha-\theta)=\frac{(\alpha-\theta)\psi'(\alpha)}{\alpha}.
\end{equation*}
This result means that the largest families at time $t$ have a size of the order of
$e^{(\alpha-\theta)t}$, and that their age is of the order of $t$ minus a constant, i.e.\ were
born in the first moments of the population growth. In particular, the largest and oldest
families are the same. 
\bigskip

This result can be interpreted as follows: If $N_\theta(t)$ denotes the size of a clonal
family started at time 0 from one individual, then conditional on its survival at time $t$,
$N_\theta(t)\, e^{-(\alpha - \theta)t}$ converges in distribution to an exponential r.v.\ with
parameter $\psi_\theta'(\alpha-\theta)$~\cite{L10}. If we restrict the limit in the last
statement to its Dirac term, we recover the previous convergence stated for the ancestral type
$$
\lim_{t\to +\infty} \PP_t (Z_0(t) > c e^{(\alpha-\theta)t}) = \frac{\alpha - \theta}{ \alpha}\ e^{-c \psi_\theta '(\alpha-\theta) }.
$$
The prefactor $\frac{\alpha - \theta}{ \alpha}$ is the probability of survival of the ancestral family conditional on the survival of the whole population, which is the ratio of $(\alpha- \theta)/b$ (survival probability of the ancestral family) with $\alpha/b$ (survival probability of the whole population).

In order to recover $\EE_t M_t(ce^{(\alpha-\theta)t},t-a_1,t-a_0)$, we need to integrate
the mutation rate per branch $\theta \, dy$ against the expected number of individuals alive
at time $y$ \emph{having at least one alive (not necessarily clonal) descendant at time $t$},
times the probability that the splitting tree spanned by one of these individuals has at least
$ce^{(\alpha-\theta)t}$ \emph{clonal} descendants, which is exactly $\PP_{t-y}(Z_0(t-y)>c
e^{(\alpha-\theta)t})$ since this splitting tree is not extinct after $t-y$ time units, which converges to 
$$
\frac{\alpha - \theta}{ \alpha}\ e^{-c \psi_\theta '(\alpha-\theta) e^{(\alpha-\theta) y}}
$$
as $t\to+\infty$. Now,
the number of individuals alive at time $y$ having descendants at time $t$ is given by the
number of branches higher than $t-y$ in the coalescent point process, i.e.\ has a geometric
distribution of parameter $\PP(H>t\mid H>t-y)=W(t-y)/W(t)$. As will appear in
Lemma~\ref{lem:W} below, this quantity converges to $e^{-\alpha y}$ as $t\rightarrow+\infty$,
which completes the interpretation of each term of~(\ref{eq:result-supercrit}).
\bigskip

Using the theory of time-dependent random characteristics, Z.\ Ta\"ib~\cite{Taib} was able to obtain
more precise results (but without considering ages of haplotypes) in the case of CMJ processes
where mutations occur \emph{at birth}. In this case, using the notation $\alpha-\theta$ for
the (positive) Malthusian parameter clonal families (to be consistent with our notation), he
obtained in Theorem~(4.6) the convergence of the number of haplotypes carried by more than $c
e^{(\alpha-\theta)t}$ individuals to a mixed Poisson r.v.\ with parameter $C w_\infty/\alpha
c^{\alpha/(\alpha-\theta)}$ where the constant $C$ is explicit and $w_\infty$ is the limit of
$N_t e^{-\alpha t}$ when $t\rightarrow+\infty$, where $N_t$ is the population size at time
$t$. This result is consistent with ours, although the precise value of the constant $C$ is
not the same.

Although not stated in~\cite{Taib}, one can easily extend this result
using~\cite[Thm.\:A.1.]{Leadbetter} to obtain the convergence in distribution on the event of
non-extinction of the point measure $\eta_t(\cdot)$, where $\eta_t([a_0,a_1])$ is the number
of haplotypes carried at time $t$ by a number of individuals in
$[a_0e^{(\alpha-\theta)t},a_1e^{(\alpha-\theta)t}]$, towards a mixed Poisson point process (also
known as \emph{Cox process}) on $\RR_+$ with intensity measure
$$
\mu(dx):=\frac{C w_\infty}{(\alpha-\theta)x^{1+\alpha/(\alpha-\theta)}}\,dx.
$$
This is actually the kind of results that we are able to prove when mutations occur during the
life of individuals and when clonal families are subcritical or critical (see below).
Unfortunately, the method we develop in Section~\ref{sec:subcrit} does not apply to the
supercritical case.

\subsubsection{Subcritical clonal families}
\label{sec:results-subcrit}

When $\alpha<\theta$, we prove in Proposition~\ref{prop:age-subcrit} that for all $a\in\RR$
$$
\lim_{t\rightarrow+\infty}\EE_t O_t\left(\frac{\alpha
    t}{\theta}+a\right)=k\,e^{-\theta a}
$$
for an explicit constant $k$, and that the maximal size of families older than $\alpha
t/\theta+a$ is tight when $t\rightarrow+\infty$ for all $a\in\RR$. We also prove in
Proposition~\ref{prop:frequent-a<t} that for all $c\in\RR$,
\begin{equation}
  \label{eq:subcrit-expect-large}
  \EE_t L_t(x_t(c))\sim k\,\varphi(\theta)^{c+\{-x_t(c)\}},  
\end{equation}
where $\{x\}$ denotes the fractional part of the real number $x$, i.e.\ $\{x\}=x-\lfloor
x\rfloor=x+\lceil-x\rceil$, where $\lfloor\cdot\rfloor$ (resp.\ $\lceil\cdot\rceil$) is the
\emph{integer part} (resp.\ \emph{ceiling}) function, and
$$
x_t(c)=k't-k''\log t+c
$$
and
\begin{equation}
  \label{eq:def-phi-theta}
  \varphi(\theta):=1-\frac{\psi(\theta)}{\theta},
\end{equation}
for explicit constants $k,k',k''$. In addition, we prove that the age of these large families
is of the order of $\log t/(\theta-\alpha)$. Hence the largest and oldest families are
different in the subcritical case.

Note that, both for large ages and large sizes, random fluctuations are of order 1 (the
parameters $a$ and $c$ are not multiplied by a function of $t$). This explains why the
right-hand side of~(\ref{eq:subcrit-expect-large}), while remaining bounded, depends on $t$:
on the one hand, the size of the largest families grows with time as $x_t(0)$ and has
fluctuations of order 1; on the other hand the size of the largest families is an integer and
hence $L_t(x_t(c))$ only depends on $\lceil x_t(c)\rceil$. Therefore, as a function of $c$, the right-hand side
of~(\ref{eq:subcrit-expect-large}) must only depend on $\lceil x_t(c)\rceil$, which is clear
since $c+\{-x_t(c)\}=-x_t(0)+\lceil x_t(c)\rceil$.

This suggests that, given any sequence of times $(t_k)_{k\geq 0}$ such that $t_k\to+\infty$ and $\{x_{t_k}(0)\}=: v$
does not depend on $k$, the r.v.\ $L_{t_k}(x_{t_k}(c))$ should converge in distribution, or
that the r.v.\ $X^{(1)}_{t_k}-x_{t_k}(0)$, where $X^{(1)}_t$ is the size of the largest family at time $t$,
should converge in distribution to some r.v.\ with values in $\ZZ-v=\{b-v, b\in\ZZ\}$. This is
what we prove in Section~\ref{sec:subcrit}. For example, we shall state here
Corollary~\ref{cor:freq-a<t}.

Let us denote by $X^{(1)}_t\geq X^{(2)}_t\geq\ldots$ the ordered sequence of family sizes in
the population at time $t$. Let also ${\cal M}(\RR)$ be the set of nonnegative $\sigma$-finite
measures on $\RR$, finite on $\RR_+$, and let us define the \emph{semi-vague} topology on
${\cal M}(\RR)$ as the one induced by all maps of the form
$$
\nu\in{\cal M}(\RR)\mapsto\int_{\RR}u(x)\nu(dx),
$$
for all bounded continuous function $u$ such that, for some $x_0\in\RR$, $u(x)=0$ for all
$x\leq x_0$. Then, Corollary~\ref{cor:freq-a<t} states that the sequence of point processes
$(Z_k)_{k\geq 0}$ on $\ZZ-v$ defined by
$$
Z_k:=\sum_{n\geq 1}\delta_{X^{(n)}_{t_k}-x_{t_k}(0)}
$$
converges in $\PP^\star$-distribution on ${\cal M}(\RR)$ equipped with the semi-vague topology
to a mixed Poisson point measure on $\ZZ-\{x_{t_0}(0)\}$ with intensity measure
$$
{\cal E}\,k\sum_{c\in \ZZ-\{x_{t_0}(0)\}}\varphi(\theta)^c\ \delta_c,
$$
where $k$ is an explicit constant and the mixture coefficient ${\cal E}$ has exponential
distribution with parameter 1.
\bigskip

We obtain similar results for the oldest families (Theorem~\ref{thm:cv-distr-old-a<t}): if
$A^{(1)}_t\geq A^{(2)}_t\geq\ldots$ denotes the ordered sequence of family ages in the
population at time $t$, the family of point processes $(Z_t,t\geq 0)$ on $\RR$ defined by
$$
Z_t:=\sum_{n\geq 1}\delta_{A_t^{(n)}-\alpha t/\theta}
$$
converges in $\PP^\star$-distribution on ${\cal M}(\RR)$ equipped with the semi-vague topology
to a mixed Poisson measure on $\RR$ with intensity measure
\begin{equation}
  \label{eq:intensity-results}
  {\cal E}\,k\,e^{-\theta a}da,  
\end{equation}
where $k$ is an explicit constant.
\bigskip

In the case of CMJ processes with mutations occuring at birth, Jagers and
Nerman~\cite[Applic.\,C]{JNb} and Taib~\cite[Prop.\,(4.2)]{Taib} obtained similar results for
the extremes of the empirical age distribution. Proposition (4.2) of~\cite{Taib} states the
convergence of the point process $Z_t$ on the event of non-extinction to a mixed Poisson
measure on $\RR$ with intensity measure
\begin{equation}
  \label{eq:intensity-results-Taib}
  w_\infty\,k'\,e^{-\theta a}da.
\end{equation}
As already noted, in the case of splitting trees, the
distribution of $w_\infty$ conditional on survival is
exponential~\cite{L10}, so that~(\ref{eq:intensity-results}) has the same form
as~(\ref{eq:intensity-results-Taib}), although the constant $k'$ is different from $k$.

The technique used in~\cite{Taib} makes use of so-called \emph{individual} time-dependent
random characteristics: for any individual $i$ in the population, we define a nonnegative
random process, called \emph{random characteristic} $(\chi_{ti}(u), u\in\RR)$, assigning some score to the individual at age $u$, such that $\chi_{ti}(u)=0$ for all $u<0$. The
random characteristic depends on an extra parameter $t$. We can then define the branching
process counted with the time-dependent random characteristic $\chi_t$ as
$$
Z^{\chi_t}_t:=\sum_{i}\chi_{ti}(t-\sigma_i),
$$
where the sum covers the set of all individuals which lived at some time in the population,
and $\sigma_i$ is the birth date of individual $i$. In such situations, the results
of~\cite{JNb} allow one to prove the convergence in distribution of $Z^{\chi_t}_t$ as
$t\rightarrow+\infty$, under a set of assumptions, among which the more stringent is that the
random characteristic is \emph{individual}, i.e.\ that for all $t\geq 0$, the random processes
$(\xi_i,\chi_{ti})$ for $i$ running in the set of all individuals in the population \emph{are
  i.i.d.}, where $(\xi_i(u),u\geq 0)$ is the process counting the number of children of
individual $i$ before age $u$.

In~\cite{Taib}, this method is applied to the population of \emph{haplotypes} (i.e.\
individuals above have to be understood as haplotypes), defining for any haplotype $i$ (a
variant of) the time-dependent random characteristic
$$
\chi_{ti}(u):=\mathbbm{1}_{\left\{\frac{\alpha
      t}{\theta}+a\leq u<\lambda_i\right\}},
$$
where $\lambda_i$ is the life length of haplotype $i$.

This method cannot be used when mutations occur during the life of individuals, since the age
of this individual when the mutation occurs influences the distribution of the progeny of the
new haplotype (except when lifetime durations are exponential r.v.), which contradicts the
assumption that $(\xi_i,\chi_{ti})$ are i.i.d. One may think of defining another random
characteristic based on individuals rather than haplotypes, counting the number of mutations
experienced by each individual, which occured more than $\alpha t/\theta+a$ time units ago and
which has descendants living at present time. With this choice, the random characteristic does
not depend on the age of the individual's mother. However, it depends on the whole progeny of
the individual, so that $(\xi_i,\chi_i)$ and $(\xi_{i'},\chi_{i'})$ are independent only if
the individuals $i$ does not descend from $i'$ and conversely. Therefore, the method
of~\cite{Taib} cannot be applied to our case.

Proposition~(4.2) of~\cite{Taib} makes use of precise estimates on the tail distribution of
the extinction time of a clonal family, which are well-known in this context. No results are
given in~\cite{Taib} on the size of large families in the subcritical case, presumably because
their method would require precise estimates on the tail distribution of the size of a clonal
population at any time, which are to our knowledge not known in general for CMJ processes. In
our case of splitting trees with mutations occurring during the life of individuals, our
formulae~(\ref{eq:exp-FS}) and~(\ref{eq:ancestral-FS}) for the expected frequency spectrum are
\emph{exact}. This allows us to obtain precise estimates for the tail distribution of the size
of a clonal population at some time, \emph{conditionally on the survival of the (clonal or not
  clonal) progeny of this population at time $t$}. We are then able to deduce exact
asymptotics for the tail distribution of the size of the largest family using a different
method than in~\cite{Taib} (see the proof of Theorem~\ref{thm:cv-loi-most-frequent}). 
\bigskip

We chose here to present results on largest and oldest haplotypes, but our method easily
applies to intermediate regimes. For example, one can easily adapt our calculations to prove
that, for all $\gamma\in[0,1]$ and $c\in\RR$,
$$
\EE_t
M_t\left(x_t(c),\frac{\alpha\gamma}{\theta}t\right)\sim
k\,\varphi(\theta)^{c+\{-x_t(c)\}},
$$
where
$$
x_t(c)=\frac{\alpha(1-\gamma)t}{-\log\varphi(\theta)}+c.
$$
Similarly, one can prove the convergence of the point process of sizes of families older than
$(\alpha\gamma/\theta)t$ as $t\rightarrow+\infty$ and compute the limit as a mixed Poisson
point measure.
\bigskip


\subsubsection{Critical clonal families}
\label{sec:results-crit}

When $\alpha=\theta$, we prove in Proposition~\ref{prop:a=t-age} that, for all $a\in\RR$,
$$
\lim_{t\rightarrow +\infty} \EE_t O_t\left(t-\frac{\log t}{\alpha}+a\right)=k e^{-\alpha
  a}
$$
for an explicit constant $k$, and that the maximal size of families older than $t-\log
t/\alpha+a$ is tight as $t\rightarrow+\infty$ for all $a\in\RR$. As in the subcritical case,
fluctuations are of order 1 here. We also prove in Proposition~\ref{prop:a=t} that, for all
$c\in\RR$,
$$
\lim_{t\rightarrow+\infty}\EE_t L_t(x_t(c))=k\,\exp\left(-\frac{2\psi'(\alpha)}{\alpha}\, c\right),
$$
where
\begin{equation}
  \label{eq:def-xt-crit-intro}
  x_t(c)=k't^2+k''t\log t+ct  
\end{equation}
and the constants $k,k',k''$ are explicit. In addition, we prove that the age of these large
families is of the order of $t/2$. Here, since the fluctuations are of the order of $t$, the
limit does not involve $\{-x_t(c)\}$ as in~(\ref{eq:subcrit-expect-large})

We are then able to deduce from these estimates the convergence of correctly rescaled point
measures of the size of the largest and the age of the oldest families
(Theorems~\ref{thm:cv-loi-freq-crit} and~\ref{thm:cv-distr-old-a=t}): using the same notation
as in Section~\ref{sec:results-subcrit}, the family of point measures $(Z_t,t\geq 0)$ defined
by
$$
Z_t:=\sum_{n\geq 1}\delta_{\frac{X^{(n)}_t-x_t(0)}{t}}
$$
converges in $\PP^\star$-distribution on ${\cal M}(\RR)$ equipped with the semi-vague topology
to a mixed Poisson measure on $\RR$ with intensity measure
$$
{\cal E}\,k\,\exp\left(-\frac{2\psi'(\alpha)}{\alpha}c\right)\,dc,
$$
where the constant $k$ is explicit.

Note that this statement and the definition~(\ref{eq:def-xt-crit-intro}) are actually a little
bit different than those of Sections~\ref{sec:a=t} and~\ref{sec:freq-crit}. However, the
results stated here can be proved by slightly modifying the proofs of
Proposition~\ref{prop:a=t} and Theorem~\ref{thm:cv-loi-freq-crit}. We  have chosen to leave this to the interested
reader.

Similarly, we obtain the convergence of the family of point measures $(Z_t,t\geq 0)$ defined
by
$$
Z_t:=\sum_{n\geq 1}\delta_{A^{(n)}_t-t+\frac{\log t}{\alpha}}
$$
to a mixed Poisson point measure on $\RR$ with intensity measure
\begin{equation}
  \label{eq:intensity-crit-intro}
  {\cal E}\,e^{-\alpha a}.  
\end{equation}

Again, in~\cite{JNb} and~\cite{Taib}, a similar result is stated only for extreme ages. It
takes the same form as our result with ${\cal E}$ replaced by $w_\infty$ and with a different
multiplicative constant in the intensity measure~(\ref{eq:intensity-crit-intro}).

\section{Large time asymptotics for the expected number of frequent or old haplotypes}
\label{sec:expected}

Our goal here is to prove the convergence results on the expectation of $L_t(x_t)$ and
$O_t(s_t)$ stated above, when clonal families are supercritical, subcritical or critical. We
start with preliminary estimates on the scale functions $W$ and $W_\theta$.

\subsection{Preliminary results on  scale functions}

\begin{lem}
  \label{lem:W}
  The survival probability of the splitting tree is $\alpha/b$, and the scale function $W$ has the following
  asymptotic behavior
  $$
  W(t)e^{-\alpha t}=\frac{1+O(e^{-\gamma t})}{\psi'(\alpha)}
  $$
  as $t\rightarrow+\infty$, for some constant $\gamma>0$.
\end{lem}

\paragraph{Proof.}
The expression of the survival probability and the fact that $W(t)\sim e^{\alpha
  t}/\psi'(\alpha)$ were already proved in~\cite{L10}. In order to get the higher-order term,
we use the fact that
\begin{equation}
  \label{eq:equiv-J}
  \PP(J>t)=\frac{\alpha}{b}+O(e^{-\gamma t})  
\end{equation}
as $t\rightarrow+\infty$, where $J$ is the extinction time of the splitting tree started from
one individual with random lifespan, distributed as $V$. Indeed, we know from \cite{L10} that
the law $\PP^\natural:=\PP(\cdot \mid J<\infty)$ is that of a subcritical splitting tree with
lifespan measure $e^{-\alpha r} \Lambda(dr)$. In particular, under $\PP^\natural$, the
lifetime $V$ of a single individual has exponential moments, and the first hitting time
$\tau_0$ of 0 by the contour process of the spliting tree also has
exponential moments (because its Laplace exponent is the inverse of $\psi^\natural (\cdot)
:=\psi(\cdot +\alpha)$). 
Since $J\le \tau_0$ a.s., $J$ has exponential moments, that is, there is some $\gamma>0$ such
that $\EE^\natural(e^{2\gamma J})<\infty$. As a consequence, also since $\alpha/b
=\PP(J=\infty)$,
$$
\PP(J>t)-\frac{\alpha}{b}= \left(1-\frac{\alpha}{b}\right)\ \PP^\natural (J>t)= O(e^{-\gamma t}).  
$$

Therefore, it follows
from~(\ref{eq:link-W-excursion}) that
$$
W(t)e^{-\alpha t}=\frac{\exp\left(-\int_t^\infty(b\PP(J>x)-\alpha)dx\right)}{\psi'(\alpha)},
$$
and the result then follows from~(\ref{eq:equiv-J}).\hfill$\Box$
\bigskip

From this result and the definition~(\ref{eq:W_theta}) of $W_\theta$, we can deduce the
following lemma. We recall that
$$
\psi_\theta(x)=\frac{x\psi(x+\theta)}{x+\theta} \qquad\text{and}\qquad
\psi'_\theta(\alpha-\theta)=\frac{(\alpha-\theta)\psi'(\alpha)}{\alpha}.
$$
\begin{lem}
  \label{lem:W_theta}
  \begin{description}
  \item[\textmd{(i)}] Assume $\alpha>\theta\geq 0$. Then
    $$
    W_\theta(t)\sim\frac{e^{(\alpha-\theta)t}}{\psi'_\theta(\alpha-\theta)}
    $$
    as $t\rightarrow+\infty$.
  \item[\textmd{(ii)}] Assume $0<\alpha<\theta$. Then
    \begin{equation}
      \label{eq:W_t-souscrit-1}
      \frac{\theta}{\psi(\theta)}-W_\theta(t)\sim
      \frac{e^{-(\theta-\alpha)t}}{|\psi'_\theta(\alpha-\theta)|}    
    \end{equation}
    as $t\rightarrow\infty$, and
    $$
    1-\frac{1}{W_\theta(t)}=\varphi(\theta)(1-\rho(t)),
    $$
    where $\varphi(\theta)$ is defined in~(\ref{eq:def-phi-theta}) and $\rho(\cdot)$ is a
    non-increasing function such that
    \begin{equation}
      \label{eq:W_t-souscrit-2}
      \rho(t)\sim\frac{\psi^2(\theta)}{\theta^2\varphi(\theta)|\psi'_\theta(\alpha-\theta)|}\,
      e^{-(\theta-\alpha)t}
    \end{equation}
    as $t\rightarrow\infty$.
  \item[\textmd{(iii)}] Assume $\alpha=\theta>0$. Then
    $$
    W_\alpha(t)=\frac{\alpha t}{\psi'(\alpha)}+\frac{1}{\psi'(\alpha)}+\alpha\int_0^{+\infty}\left(W(y)e^{-\alpha
        y}-\frac{1}{\psi'(\alpha)}\right)dy+o(1)
    $$
    as $t\rightarrow+\infty$.
  \end{description}
\end{lem}

\paragraph{Proof.}
Points~(i) and~(iii) are easy consequences of Lemma~\ref{lem:W} and the
definition~(\ref{eq:W_theta}) of $W_\theta$. Point~(i) can also be seen as a trivial corollary
of Lemma~\ref{lem:W} since when $\alpha>\theta$, $W_\theta$ is the scale function of a
supercritical splitting tree.

For Point~(ii), by Tauberian theorems (see~\cite{L10}), we have $W_\theta(t)\rightarrow
1/\psi'_\theta(0)=\theta/\psi(\theta)$ as $t\rightarrow\infty$. Since
$$
\frac{1}{\psi(\theta)}=\int_0^\infty W(y)e^{-\theta y}dy ,
$$
(\ref{eq:W_theta}) yields
$$
W_\theta(t) - \frac{\theta}{\psi(\theta)} = e^{-\theta t} W(t) - \theta \int_t^\infty W(y)e^{-\theta y}dy .
$$
Since $W(t)\sim e^{\alpha t}/\psi'(\alpha)$, we have
$$
\int_t^\infty W(y)e^{-\theta y}dy \,\sim\, \frac{1}{\psi'(\alpha)}\int_t^\infty e^{-(\theta-\alpha)y}dy
$$
as $t\rightarrow+\infty$. This entails~(\ref{eq:W_t-souscrit-1}). Since one has
$$
\rho(t)=\frac{1}{\varphi(\theta)}\left(\frac{1}{W_\theta(t)}-\frac{\psi(\theta)}{\theta}\right)
$$
and
$$
\frac{1}{W_\theta(t)}-\frac{\psi(\theta)}{\theta}
=\frac{\frac{\theta}{\psi(\theta)}-W_\theta(t)}{W_\theta(t)\,\frac{\theta}{\psi(\theta)}},
$$
the proof of~(ii) is easily completed.\hfill$\Box$ 
\bigskip

\subsection{The case of supercritical families ($\alpha>\theta$)}

In the case of supercritical clonal families, the asymptotic expected number of frequent
haplotypes can be explicitly computed. Note that in the next statement and elsewhere in the
paper, the Dirac measure at time 0 corresponds to the contribution of the family carrying the
ancestral type.

\begin{prop}
  \label{prop:frequent-a>t}
  Assume $\alpha>\theta\geq 0$ and let $0\leq a_0<a_1\leq+\infty$ and $c\geq 0$. For all $t\geq 0$, let
  $$
  x_t(c)=c\exp((\alpha-\theta)t).
  $$
  Then
  $$
  \lim_{t\rightarrow+\infty}\EE_t[M_t(x_t(c),t-a_1,t-a_0)]=\frac{\alpha-\theta}{\alpha}\int_{a_0}^{a_1}
  \exp\left(\alpha y
    -c\,\psi'_\theta(\alpha-\theta)\,e^{(\alpha-\theta)y}\right)(\theta\:dy+\delta_0(dy)).
  $$
\end{prop}

\paragraph{Proof.}
Using~(\ref{eq:ancestral-FS}),~(\ref{eq:exp-FS-with-age}) and Lemma~\ref{lem:W}, for all
$t\geq a$, we have
\begin{align}
  \EE_t & [M_t(x_t,t-b,t-a)]=\sum_{k=\lceil
    ce^{(\alpha-\theta)t}\rceil}^{+\infty}\left[
   W(t)\frac{e^{-\theta t}}{W_\theta(t)^2}\left(1-\frac{1}{W_\theta(t)}\right)^{k-1}\mathbbm{1}_{\{a=0\}}
   \right. \notag \\ &
  \qquad\qquad\qquad\qquad\qquad\qquad\qquad\qquad\qquad\qquad
  \left.+
   \int_{(t-b)\vee 0}^{t-a}\theta
    W(t)\frac{e^{-\theta x}}{W_\theta(x)^2}\left(1-\frac{1}{W_\theta(x)}\right)^{k-1}dx
    \right]
     \notag \\
  & =\mathbbm{1}_{\{a=0\}}\frac{W(t)}{W_\theta(t)}e^{-\theta
    t}\left(1-\frac{1}{W_\theta(t)}\right)^{\lceil ce^{(\alpha-\theta)t}\rceil-1} 
  +W(t)\int_{(t-b)\vee 0}^{t-a}\theta\frac{e^{-\theta x}}{W_\theta(x)}\left(1-\frac{1}{W_\theta(x)}\right)^{\lceil
    ce^{(\alpha-\theta)t}\rceil-1}\:dx \notag \\
  & \sim\frac{1}{\psi'(\alpha)}\int_{a}^{b\wedge t}\frac{e^{\theta y}}{W_\theta(t-y)e^{-(\alpha-\theta)t}}
  \exp\left((\lceil
    ce^{(\alpha-\theta)t}\rceil-1)\log(1-1/W_\theta(t-y))\right)(\theta dy+\delta_0(dy)),
  \label{eq:calcul-esp+gde-famille}
\end{align}
where we recall that $\lceil\cdot\rceil$ is the ceiling function. Since $W_\theta$ is
nondecreasing and $W_\theta(0)=1$, it follows from Lemma~\ref{lem:W_theta}~(i) that there
exists a constant $C>0$ such that
$$
\frac{1}{C}e^{(\alpha-\theta)t}\leq W_\theta(t)\leq C
e^{(\alpha-\theta)t},\qquad\forall t\geq 0.
$$
Therefore, for all $y\geq 0$, the quantity inside the integral in the
r.h.s.\ of~(\ref{eq:calcul-esp+gde-famille}) is smaller than
$$
C'e^{\alpha y}\exp\left(-C'' e^{(\alpha-\theta)y}\right)
$$
for some constants $C',C''>0$. Now, using Lemma~\ref{lem:W_theta}~(i)
again, for all $y\geq 0$, the quantity inside the integral in the
r.h.s.\ of~(\ref{eq:calcul-esp+gde-famille}) converges to
$$
\frac{\psi'(\alpha)(\alpha-\theta)}{\alpha}e^{\alpha
  y}\exp\left(-c \psi'_\theta(\alpha-\theta) e^{(\alpha-\theta)y}\right)
$$
when $t\rightarrow+\infty$.  Proposition~\ref{prop:frequent-a>t} then
follows from the dominated convergence theorem.\hfill$\Box$
\bigskip

\subsection{The case of subcritical families ($\alpha<\theta$)}

Our first result deals with ages and sizes of the oldest clonal families. Note that the
scaling constant $a$ varies in $\RR$.
\begin{prop}
  \label{prop:age-subcrit} When $\alpha<\theta$, for any $a\in\RR$,
  \begin{equation}
    \label{eq:old-subcrit}
    \lim_{t\rightarrow+\infty}\EE_t\left[O_t\left(\frac{\alpha\,t}{\theta}+a\right)\right]
    =\frac{\psi(\theta)}{\theta\psi'(\alpha)}e^{-\theta a}.
  \end{equation}
  In addition, for any $x_t\rightarrow+\infty$ as $t\rightarrow+\infty$,
  \begin{equation}
    \label{eq:old-subcrit-2}
  \lim_{t\rightarrow+\infty}\EE_t\left[M_t\left(x_t,\frac{\alpha\,t}{\theta}+a\right)\right]=0.
  \end{equation}
\end{prop}

\paragraph{Proof.}
Using~(\ref{eq:ancestral-FS}) and~(\ref{eq:exp-FS-with-age}) as in the proof of
Proposition~\ref{prop:frequent-a>t}, we have
\begin{align*}
  \EE_t[O_t(\alpha t/\theta+a)] & = W(t)\int_{\frac{\alpha\,t}{\theta}+a}^t\frac{e^{-\theta
      x}}{W_\theta(x)}(\theta dx+\delta_t(dx)) \\ & \sim\frac{\psi(\theta)e^{\alpha
      t}}{\psi'(\alpha)}\int_{\frac{\alpha\,t}{\theta}+a}^t e^{-\theta x}\Big(dx+\frac{1}{\theta}\delta_t(dx)\Big) ,
\end{align*}
where we used that $W_\theta(x)\rightarrow\theta/\psi(\theta)$ as
$x\rightarrow+\infty$ (Lemma~\ref{lem:W_theta}). Eq.~(\ref{eq:old-subcrit}) then
easily follows.

Similarly,
\begin{align*}
  \EE_t[M_t(x_t,\,\alpha t/\theta+a)] & = W(t)\int_{\frac{\alpha\,t}{\theta}+a}^t\frac{e^{-\theta
      x}}{W_\theta(x)}\left(1-\frac{1}{W_\theta(x)}\right)^{\lceil
    x_t\rceil-1}(\theta dx+\delta_t(dx)) \\ & 
  \leq \left(1-\frac{1}{W_\theta(t)}\right)^{\lceil
    x_t\rceil-1}
  W(t)\int_{\frac{\alpha\,t}{\theta}+a}^t\frac{e^{-\theta
      x}}{W_\theta(x)}(\theta dx+\delta_t(dx)),
\end{align*}
since $W_\theta$ is nondecreasing. Eq.~(\ref{eq:old-subcrit-2}) then follows
from~(\ref{eq:old-subcrit}) and the fact that $(1-1/W_\theta(t))^{\lceil x_t\rceil}\rightarrow
0$ as $t\rightarrow+\infty$, since $W_\theta(t)\rightarrow \theta/\psi(\theta)>1$.\hfill$\Box$
\bigskip

Our next result gives the asymptotics of the expected number of large families (notice again
that the scaling constant $c$ varies in $\RR$) and states that their ages all are
asymptotically equivalent to $(\theta-\alpha)^{-1}\log(t)$.
We recall that $\varphi(\theta)<1$, so that
$|\log\varphi(\theta)|=-\log\varphi(\theta)$.

\begin{prop}
  \label{prop:frequent-a<t}
  Assume $\alpha<\theta$. For all $c\in\RR$, let
  $$
  x_t(c):=\frac{\alpha t- \frac{\theta}{\theta-\alpha}\log t}{|\log\varphi(\theta)|}+c.
  $$
  Then, for all $\varepsilon>0$,
  \begin{equation}
    \label{eq:frequent-a<t}
    \EE_t[L_t(x_t(c))] \sim \EE_t\left[M_t\left(x_t(c),\, \frac{1-\varepsilon}{\theta-\alpha}\log
        t,\,\frac{1+\varepsilon}{\theta-\alpha}\log t\right)\right] \sim
    A(\theta)\,\varphi(\theta)^{c-1+\{-x_t(c)\}}
  \end{equation}
  as $t\rightarrow+\infty$, where we recall that $\{\cdot\}$ is the fractional part function
  and
  \begin{equation}
    \label{eq:def-A}
    A(\theta):=\frac{\Gamma\left(\frac{\theta}{\theta-\alpha}\right)\psi(\theta)}{\alpha}\,
    |\psi'_\theta(\alpha-\theta)|^{\frac{\alpha}{\theta-\alpha}}
    \left(\frac{\theta^2}{\alpha \psi(\theta)^2}\,\varphi(\theta)\ |\log\varphi(\theta)|
      \right)^{\frac{\theta}{\theta-\alpha}},
  \end{equation}
  where $\Gamma(s)=\int_0^\infty y^{s-1}e^{-y}dy$ is the Gamma function.
\end{prop}

\paragraph{Proof.}
Proceeding similarly as in the proof of
Proposition~\ref{prop:frequent-a>t}, we deduce from Lemma~\ref{lem:W_theta}~(ii) that
\begin{align}
  \EE_t[L_t(x_t(c))] & =W(t)\,\varphi(\theta)^{\lceil x_t(c)\rceil-1}
  \int_0^t\frac{e^{-\theta x}}{W_\theta(x)}\,(1-\rho(x))^{\lceil x_t(c)\rceil-1}\, (\theta dx+\delta_t(dx))
  \notag \\ & \sim\frac{t^{\theta/(\theta-\alpha)}\,\varphi(\theta)^{c-1+\{-x_t(c)\}}}{\psi'(\alpha)}
  \int_0^t\frac{e^{-\theta x+(\lceil x_t(c)\rceil-1)\log(1-\rho(x))}}{W_\theta(x)} (\theta
  dx+\delta_t(dx)) \label{eq:calcul-pas-drole}
\end{align}
as $t\rightarrow+\infty$, where we used the relation $\lceil x_t(c)\rceil-1=x_t(0)+c-1+\{-x_t(c)\}$. We will now abridge $x_t(c)$ into $x_t$.  

Let $\varepsilon >0$. Let us first bound from above the previous integral restricted to the
complement of $\left[\frac{1-\varepsilon}{\theta-\alpha}\log t,
  \frac{1+\varepsilon}{\theta-\alpha}\log t\right]$.
On the one hand, using the inequality
$\log(1-\rho(x))\leq 0$ and the fact that $W_\theta(x)$ converges to a positive constant when
$x\rightarrow+\infty$, we have for all $\varepsilon>0$
$$
\int_{\frac{1+\varepsilon}{\theta-\alpha}\log t}^{+\infty}\frac{e^{-\theta x+(\lceil
    x_t\rceil-1)\log(1-\rho(x))}}{W_\theta(x)} (\theta dx+\delta_t(dx))\leq C\left(e^{-\theta
  t}+e^{-\theta\frac{1+\varepsilon}{\theta-\alpha}\log t}\right)=o(t^{-\frac{\theta}{\theta-\alpha}})
$$
as $t\rightarrow+\infty$.

On the other hand, since $\rho(x)$ is non-increasing and using~(\ref{eq:W_t-souscrit-2}),
$\rho(x)\geq\rho(\frac{1-\varepsilon}{\theta-\alpha}\log t)\geq Ct^{-1+\varepsilon}$ for all $0\leq x\leq
(1-\varepsilon)/(\theta-\alpha)\:\log t$. Therefore,
$$
\int_{0}^{\frac{1-\varepsilon}{\theta-\alpha}\log t}\frac{e^{-\theta x+(\lceil
    x_t\rceil-1)\log(1-\rho(x))}}{W_\theta(x)} (\theta dx+\delta_t(dx))\leq
C'\int_0^{\frac{1-\varepsilon}{\theta-\alpha}\log t}e^{-\theta x}\,e^{-C(\lceil x_t\rceil-1)
  t^{-1+\varepsilon}}dx.
$$
Since $(\lceil x_t\rceil-1) t^{-1+\varepsilon}\geq C''t^\varepsilon$ for $t$ large enough, we
deduce that the previous integral is also $o(t^{-\frac{\theta}{\theta-\alpha}})$. In
conclusion, $\EE_t[L_t(x_t(a))]$ and $\EE_t\left[M_t\left(x_t(a),\,
    \frac{1-\varepsilon}{\theta-\alpha}\log t,\, \frac{1+\varepsilon}{\theta-\alpha}\log
    t\right)\right]$ are both asymptotically equivalent to $E_t(\varepsilon)$, provided that
$E_t(\varepsilon)$ is uniformly  bounded from below, where
\begin{equation}
  \label{eq:goal}
  E_t(\varepsilon):=
  \frac{\psi(\theta)\,t^{\theta/(\theta-\alpha)}\,\varphi(\theta)^{a-1+\{-x_t(a)\}}}{\psi'(\alpha)}
  \int_{\frac{1-\varepsilon}{\theta-\alpha}\log t}^{\frac{1+\varepsilon}{\theta-\alpha}\log t}
  e^{-\theta x+(\lceil x_t\rceil-1)\log(1-\rho(x))} dx.
\end{equation}
Note that, since $W_\theta(x)\rightarrow\psi(\theta)/\theta$ when $x\rightarrow+\infty$, the
replacement of $W_\theta(x)$ by its limit in the r.h.s.\ of~(\ref{eq:calcul-pas-drole}) is
justified. The proof of Proposition~\ref{prop:frequent-a<t} will hence be completed if we can
prove that
\begin{equation*}
  E_t(\varepsilon)\sim A(\theta)
  \,\varphi(\theta)^{c-1+\{-x_t(c)\}}.
\end{equation*}
This is the aim of the rest of the proof. Set
\begin{equation}
  \label{eq:def-B-pf-sscrit}
  B:=\frac{\psi^2(\theta)}{\theta^2\varphi(\theta)|\psi'_\theta(\alpha-\theta)|},  
\end{equation}
and we make the change of variable $y:=B(\lceil x_t\rceil-1) e^{-(\theta-\alpha)x}$ in~(\ref{eq:goal}):
\begin{multline}
  E_t(\varepsilon)\sim
  \frac{\psi(\theta)\,t^{\theta/(\theta-\alpha)}\,\varphi(\theta)^{c-1+\{-x_t\}}}{\psi'(\alpha)(\theta-\alpha)(B(\lceil
  x_t\rceil-1))^{\theta/(\theta-\alpha)}}
  \int_{B(\lceil x_t\rceil-1) t^{-1-\varepsilon}}^{B(\lceil x_t\rceil-1) t^{-1+\varepsilon}}
  y^{\theta/(\theta-\alpha)-1} \\
  \exp\biggl((\lceil x_t\rceil-1)\log\biggl(1-\rho\biggl(\frac{\log\frac{B(\lceil
      x_t\rceil-1)}{y}}{\theta-\alpha}\biggr)\biggr)\biggr)\,dy.
  \label{eq:calcul-bourrin}
\end{multline}
Now, using Lemma~\ref{lem:W_theta}~(ii) again, we have for all $y>0$
$$
\rho\biggl(\frac{\log\frac{B(\lceil
      x_t\rceil-1)}{y}}{\theta-\alpha}\biggr)\sim Be^{\log\frac{y}{B(\lceil x_t\rceil-1)}}=\frac{y}{\lceil
    x_t\rceil-1}
$$
as $t\rightarrow+\infty$. Therefore, the exponential in the r.h.s.\ of~(\ref{eq:calcul-bourrin}) converges
for all $y>0$ to $e^{-y}$ when $t\rightarrow+\infty$. In addition, using the inequalities
$\log(1-\rho(x))\leq-\rho(x)$ and $\rho(x)\geq Ce^{-(\theta-\alpha)x}$ for all $x$ large enough, we have
$$
y^{\frac{\theta}{\theta-\alpha}-1}\exp\biggl((\lceil
x_t\rceil-1)\log\biggl(1-\rho\biggl(\frac{\log\frac{B(\lceil
      x_t\rceil-1)}{y}}{\theta-\alpha}\biggr)\biggr)\biggr)\leq y^{\frac{\alpha}{\theta-\alpha}}e^{-Cy}.
$$
Since $x_t t^{-1-\varepsilon}\rightarrow 0$ and $x_t t^{-1+\varepsilon}\rightarrow+\infty$ when
$t\rightarrow+\infty$, Lebesgue's theorem finally yields
$$
E_t(\varepsilon)\sim \frac{\psi(\theta)\varphi(\theta)^{c-1+\{-x_t\}}
  \Gamma(\frac{\theta}{\theta-\alpha})} {\psi'(\alpha)(\theta-\alpha)B^{\theta/(\theta-\alpha)}}\,
\left(\frac{t}{\lceil x_t\rceil-1}\right)^{\theta/(\theta-\alpha)} ,
$$
Remembering that $x_t\sim \alpha t/|\log\varphi(\theta)|$ as $t\rightarrow+\infty$
concludes the proof of Proposition~\ref{prop:frequent-a<t}.\hfill$\Box$
\bigskip

\subsection{The case of critical families ($\alpha=\theta$)}
\label{sec:a=t}

The following result gives the asymptotic expected number of old families and states that
their sizes are tight.

\begin{prop}
  \label{prop:a=t-age} Assume $\alpha=\theta>0$. For all $a\in\RR$, we have
  $$
  \lim_{t\rightarrow+\infty}\EE_t\left[O_t\left(t-\frac{\log t}{\alpha}+a\right)\right]=\frac{e^{-\alpha
      a}}{\alpha}.
  $$
  In addition, for all $x_t\rightarrow+\infty$,
  $$
  \lim_{t\rightarrow+\infty}\EE_t\left[M_t\left(x_t,\,t-\frac{\log t}{\alpha}+a\right)\right]=0.
  $$
\end{prop}

\paragraph{Proof.}
Using Lemma~\ref{lem:W_theta}~(iii), a similar computation as for Proposition~\ref{prop:age-subcrit} yields
$$
\EE_t[O_t(t-\log t/\alpha+a)]\sim e^{\alpha t}\int_{t-\frac{\log t}{\alpha}+a}^t\frac{e^{-\alpha
    x}}{x}(dx+\frac{1}{\alpha}\delta_t(dx))\sim \frac{e^{\alpha t}}{t}\int_{t-\frac{\log
    t}{\alpha}+a}^te^{-\alpha x}dx+\frac{1}{\alpha t}.
$$
The first limit easily follows. The second limit is obtained exactly as in the proof of
Proposition~\ref{prop:age-subcrit}.\hfill$\Box$ 
\bigskip

The computations for the most frequent haplotypes are more involved. The following result
gives the asymptotics of the expected number of large families and states that their ages are
all asymptotically equivalent to $t/2$.

\begin{prop}
  \label{prop:a=t} Assume $\alpha=\theta>0$. For all $c\in\RR$, we define
  $$
  x_t(c):=\frac{\alpha^2}{4\psi'(\alpha)}\left(t-\frac{\log t}{2\alpha}+c\right)^2.
  $$
  Then, for all $\varepsilon>0$,
  $$
  \lim_{t\rightarrow+\infty}\EE_t[L_t(x_t(c))]=\lim_{t\rightarrow+\infty}
  \EE_t\left[M_t\left(x_t(c),\,\frac{1-\varepsilon}{2}\,t,\,
      \frac{1+\varepsilon}{2}\,t\right)\right]=\sqrt{\frac{2\pi}{\alpha}}\,e^{B-\frac{\psi'(\alpha)}{2}}\,
  e^{-\alpha\, c},
  $$
  where
  \begin{equation}
    \label{eq:def-B-crit}
    B=1+\alpha\int_0^{+\infty}\left(\psi'(\alpha)W(y)e^{-\alpha y}-1\right)dy.    
  \end{equation}
\end{prop}


\paragraph{Proof.}
Similarly as in the proof of Proposition~\ref{prop:frequent-a<t}, we have
\begin{equation}
  \label{eq:pf-crit-lim-esp}
  \EE_t[L_t(x_t(c))]\sim\frac{e^{\alpha t}}{\psi'(\alpha)}\int_0^t\frac{e^{-\alpha x}}{W_\alpha(x)}
  e^{(\lceil x_t(c)\rceil-1)\log\left(1-\frac{1}{W_\alpha(x)}\right)}(\alpha dx+\delta_t(dx)).
\end{equation}

Let $\varepsilon>0$. Let us first bound from above the previous integral restricted to the
complement of $\left[\frac{1-\varepsilon}{2} t, \frac{1+\varepsilon}{2} t\right]$.
Fix $\eta\in(0,1)$. By Lemma~\ref{lem:W_theta}~(iii), for all $x$ large enough, $1\leq
W_\alpha(x)\leq\frac{\alpha t}{\psi'(\alpha)(1-\eta)}$. Hence, using the fact that $x_t(c)\sim
\alpha^2 t^2/4\psi'(\alpha)$, for $t$ sufficiently large,
\begin{align*}
  \int_{\frac{1+\varepsilon}{2}t}^t\frac{e^{-\alpha x}}{W_\alpha(x)} e^{(\lceil
    x_t(c)\rceil-1)\log\left(1-\frac{1}{W_\alpha(x)}\right)}(\alpha dx+\delta_t(dx)) & \leq
  \int_{\frac{1+\varepsilon}{2}t}^t e^{-\alpha x-(\lceil
    x_t(c)\rceil-1)\frac{(1-\eta)\psi'(\alpha)}{\alpha x}}(\alpha dx+\delta_t(dx)) \\ & \leq
  \alpha\int_{\frac{1+\varepsilon}{2}t}^te^{-\alpha\left(x+\frac{(1-2\eta)t^2}{4x}\right)}dx
  +e^{-\alpha t(1+(1-2\eta)/4)}.
\end{align*}
The quantity inside the integral in the integral of the r.h.s.\ is maximal for
$x=(1-2\eta)^{1/2}t/2$, which is outside the integration domain, so that
$$
\int_{\frac{1+\varepsilon}{2}t}^te^{-\alpha\left(x+\frac{(1-2\eta)t^2}{4x}\right)}dx\leq
t\exp\left(-\alpha\left(\frac{t(1+\varepsilon)}{2}+\frac{t(1-2\eta)}{2(1+\varepsilon)}\right)\right).
$$
This last quantity is $o(e^{-\alpha t})$ if one chooses $\eta<\varepsilon^2/2$. 

Using the fact that $W_\alpha(x)$ is non-decreasing, larger than 1 and that $e^{-\alpha x}\leq
1$, we have
$$
\int_{0}^{\log t}\frac{e^{-\alpha x}}{W_\alpha(x)} e^{(\lceil
  x_t(c)\rceil-1)\log\left(1-\frac{1}{W_\alpha(x)}\right)}\alpha dx\leq \alpha e^{-\frac{\lceil
    x_t(c)\rceil-1}{W_\alpha(\log t)}}\log t.
$$
Since $x_t(c)\sim \alpha^2 t^2/4\psi'(\alpha)$, using Lemma~\ref{lem:W_theta}~(iii) again, for
$t$ large enough, the previous integral is smaller than
$$
\alpha\log ( t) \, e^{-Ct^2/\log t}=o(e^{-\alpha t}),
$$
for a constant $C>0$ independent of $t$.

Finally, using Lemma~\ref{lem:W_theta}~(iii) similarly as above, for all $\eta\in(0,1)$, for
$t$ large enough,
$$
\int_{\log t}^{\frac{1-\varepsilon}{2}t}\frac{e^{-\alpha x}}{W_\alpha(x)} e^{(\lceil
  x_t(c)\rceil-1)\log\left(1-\frac{1}{W_\alpha(x)}\right)}\alpha \,dx\leq
\alpha\int_{\log t}^{\frac{1-\varepsilon}{2}t}e^{-\alpha\left(x+\frac{(1-2\eta)t^2}{4x}\right)}dx.
$$
Taking $\eta$ small enough so that $(1-\eta)^{1/2}t/2>t(1-\varepsilon)/2$, the previous
quantity can be bounded from above by
$$
\alpha te^{-\frac{\alpha t}{2}\left(1-\varepsilon+\frac{1-2\eta}{1-\varepsilon}\right)}=o(e^{-\alpha t})
$$
if one takes again $\eta<\varepsilon^2/2$. In conclusion, $\EE_t[L_t(x_t(c))]$ and
$\EE_t\left[M_t\left(x_t(c),\, \frac{1-\varepsilon}{2}t,\, \frac{1+\varepsilon}{2}
    t\right)\right]$ are both asymptotically equivalent to $E_t (\varepsilon)$, provided that
$E_t(\varepsilon)$ is uniformly bounded from below, where
$$
E_t(\varepsilon):=\frac{\alpha e^{\alpha t}}{\psi'(\alpha)}\int_{\frac{1-\varepsilon}{2}t}^{\frac{1+\varepsilon}{2}t}\frac{e^{-\alpha x}}{W_\alpha(x)}
  e^{(\lceil x_t(c)\rceil-1)\log\left(1-\frac{1}{W_\alpha(x)}\right)}\,dx.
$$
Therefore, it only remains to prove that
$$
\lim_{t\rightarrow+\infty}E_t(\varepsilon)=\sqrt{\frac{2\pi}{\alpha}}\,e^{B-\frac{\psi'(\alpha)}{2}}\,
e^{-\alpha\, c}.
$$
Using the facts that $W_\alpha(x)\sim \alpha x/\psi'(\alpha)$ and that
$|\log(1-1/W_\alpha(x))|\leq C/t$ for all $x$ large enough, we have
\begin{equation*}
  E_t(\varepsilon)
  \sim e^{\alpha
    t}\int_{\frac{1-\varepsilon}{2}t}^{\frac{1+\varepsilon}{2}t} e^{-\alpha
    x}\,e^{x_t(c)\log\left(1-\frac{1}{W_\alpha(x)}\right)}\,\frac{dx}{x},
\end{equation*}
when $t\rightarrow+\infty$.

It follows from Lemma~\ref{lem:W_theta}~(iii) that
\begin{equation}
  \label{eq:dl-log-1-1/W_a}
  \log\left(1-\frac{1}{W_\alpha(x)}\right)=-\frac{\psi'(\alpha)}{\alpha
    x}+\frac{\psi'(\alpha)(B-\psi'(\alpha)/2)}{\alpha^2 x^2}+o\left(\frac{1}{x^2}\right)
\end{equation}
as $x\rightarrow+\infty$. Therefore,
$$
F(t):=x_t(c)\ \sup_{\frac{1-\varepsilon}{2}t\leq x\leq
  \frac{1+\varepsilon}{2}t}\ \left[\log\left(1-\frac{1}{W_\alpha(x)}\right)+\frac{\psi'(\alpha)}{\alpha
    x}-\frac{\psi'(\alpha)(B-\psi'(\alpha)/2)}{\alpha^2 x^2}\right]\rightarrow 0
$$
as $t\rightarrow+\infty$. Hence, using the facts that $x_t(c)\sim\frac{\alpha^2
  t}{4\psi'(\alpha)}$ and that $x\in[\frac{1-\varepsilon}{2}t,\frac{1+\varepsilon}{2}t]$, for
$t$ large enough,
\begin{multline}
  \frac{2e^{-C\varepsilon}}{(1+\varepsilon)t}\,
  e^{B-\frac{\psi'(\alpha)}{2}}\,e^{\alpha t}\int_{\frac{1-\varepsilon}{2}t}^{\frac{1+\varepsilon}{2}t}
  e^{-\alpha x-\frac{x_t(c)\psi'(\alpha)}{\alpha x}}dx \leq E_t(\varepsilon) \\ \leq
  \frac{2e^{C\varepsilon}}{(1-\varepsilon)t}\,
  e^{B-\frac{\psi'(\alpha)}{2}}\,e^{\alpha t}\int_{\frac{1-\varepsilon}{2}t}^{\frac{1+\varepsilon}{2}t}
  e^{-\alpha x-\frac{x_t(c)\psi'(\alpha)}{\alpha x}}dx, \label{eq:calcul-bourrin-2}
\end{multline}
for a constant $C$ independent of $\varepsilon$ and $t$.

Now, let us compute the asymptotic behavior of the integral involved in these inequalities:
first, the change of variable $x=\beta_t y$ with $\beta_t=\sqrt{x_t(c)\psi'(\alpha)}/\alpha$
yields
$$
\int_{\frac{1-\varepsilon}{2}t}^{\frac{1+\varepsilon}{2}t}
  e^{-\alpha x-\frac{x_t(c)\psi'(\alpha)}{\alpha x}}dx
  =\beta_t\int_{\frac{1-\varepsilon}{2\beta_t}t}^{\frac{1+\varepsilon}{2\beta_t}t}
  e^{-\alpha \beta_t\left(y+\frac{1}{y}\right)}dy.
$$
Next, we introduce the new change of variable
$$
y=\frac{2+z^2+z\sqrt{z^2+4}}{2}.
$$
This defines a $C^1$-diffeomorphism from $z\in(-\infty,+\infty)$ to $y\in(0,+\infty)$ such that
$$
y+\frac{1}{y}=2+z^2.
$$
Note that $z>0$ if and only if $y>1$, which means that
$$
z=\text{sgn}(y-1)\sqrt{y+\frac{1}{y}-2},
$$
where $\text{sgn}(x)=1$ if $x\geq 0$ and $-1$ if $x<0$. Since 
\begin{equation}
  \label{eq:def-beta-t}
  \beta_t=\frac{t-\frac{\log t}{2\alpha}+c}{2}\sim\frac{t}{2},
\end{equation}
the inequality $(1-\varepsilon)t/2\beta_t<1<(1+\varepsilon)t/2\beta_t$ holds for $t$ large enough,
which yields
$$
\int_{\frac{1-\varepsilon}{2}t}^{\frac{1+\varepsilon}{2}t}
  e^{-\alpha x-\frac{x_t(c)\psi'(\alpha)}{\alpha x}}dx
  =\beta_t\int_{-\sqrt{\frac{(1-\varepsilon)t}{2\beta_t}
      +\frac{2\beta_t}{(1-\varepsilon)t}-2}}^{\sqrt{\frac{(1+\varepsilon)t}{2\beta_t}
      +\frac{2\beta_t}{(1+\varepsilon)t}-2}}
  e^{-\alpha\beta_t(2+z^2)}\left(z+\frac{z^2+2}{\sqrt{z^2+4}}\right)\,dz.
$$
Now,
$$
\lim_{t\rightarrow+\infty}\sqrt{\frac{(1-\varepsilon)t}{2\beta_t}
      +\frac{2\beta_t}{(1-\varepsilon)t}-2}=\frac{\varepsilon}{\sqrt{1-\varepsilon}}\qquad \text{and}\qquad
\lim_{t\rightarrow+\infty}\sqrt{\frac{(1+\varepsilon)t}{2\beta_t}
      +\frac{2\beta_t}{(1+\varepsilon)t}-2}=\frac{\varepsilon}{\sqrt{1+\varepsilon}}.
$$
Since $z+\frac{z^2+2}{\sqrt{z^2+4}}$ is $C^1$ in the neighborhood of 0, with value $1$ at $z=0$, we obtain
\begin{multline*}
  (1-C\varepsilon)\beta_t e^{-2\alpha\beta_t}\int_{-\sqrt{\frac{(1-\varepsilon)t}{2\beta_t}
      +\frac{2\beta_t}{(1-\varepsilon)t}-2}}^{\sqrt{\frac{(1+\varepsilon)t}{2\beta_t}
      +\frac{2\beta_t}{(1+\varepsilon)t}-2}}
  e^{-\alpha\beta_t z^2}dz\leq\int_{\frac{1-\varepsilon}{2}t}^{\frac{1+\varepsilon}{2}t}
  e^{-\alpha x-\frac{x_t(c)\psi'(\alpha)}{\alpha x}}dx \\ \leq
  (1+C\varepsilon)\beta_t e^{-2\alpha\beta_t}\int_{-\sqrt{\frac{(1-\varepsilon)t}{2\beta_t}
      +\frac{2\beta_t}{(1-\varepsilon)t}-2}}^{\sqrt{\frac{(1+\varepsilon)t}{2\beta_t}
      +\frac{2\beta_t}{(1+\varepsilon)t}-2}}
  e^{-\alpha\beta_t z^2}dz
\end{multline*}
for $t$ large enough. Making the last change of variable $u=\sqrt{2\alpha\beta_t}\,z$ finally yields
$$
(1-C'\varepsilon)\sqrt{\frac{\pi\beta_t}{\alpha}}e^{-2\alpha\beta_t}\leq
\int_{\frac{1-\varepsilon}{2}t}^{\frac{1+\varepsilon}{2}t}
  e^{-\alpha x-\frac{x_t(c)\psi'(\alpha)}{\alpha x}}dx\leq
(1+C'\varepsilon)\sqrt{\frac{\pi\beta_t}{\alpha}}e^{-2\alpha\beta_t}
$$
for $t$ large enough. Combining this with~(\ref{eq:calcul-bourrin-2}), we obtain that, for all
$\varepsilon>0$ small enough, there exists $t_0>0$ such that, for all $t>t_0$,
$$
(1-C''\varepsilon)\frac{2}{t}e^{B-\frac{\psi'(\alpha)}{2}}e^{\alpha
  t}\sqrt{\frac{\pi\beta_t}{\alpha}}e^{-2\alpha\beta_t}\leq E_t(\varepsilon)\leq 
(1+C''\varepsilon)\frac{2}{t}e^{B-\frac{\psi'(\alpha)}{2}}e^{\alpha
  t}\sqrt{\frac{\pi\beta_t}{\alpha}}e^{-2\alpha\beta_t},
$$
where the constant $C''$ is independent of $\varepsilon$ and $t$. It then follows from~(\ref{eq:def-beta-t})
that
$$
(1-C''\varepsilon)\sqrt{\frac{2\pi}{\alpha}}e^{B-\frac{\psi'(\alpha)}{2}}\,e^{-\alpha c}\leq E_t(\varepsilon)\leq 
(1+C''\varepsilon)\sqrt{\frac{2\pi}{\alpha}}e^{B-\frac{\psi'(\alpha)}{2}}\,e^{-\alpha c}.
$$

Since we have shown in the beginning of the proof that $E_t(\varepsilon)=E_t(\varepsilon')+o(1)$ for all
$\varepsilon'<\varepsilon$, the previous inequality applied to $E_t(\varepsilon')$ concludes the proof of
Proposition~\ref{prop:a=t}.\hfill$\Box$
\bigskip

\section{Large or old families: convergence in distribution for subcritical clonal families}
\label{sec:subcrit}


In all this section, we assume $\alpha<\theta$. Our goal is to compute the joint limiting
distribution when $t\rightarrow+\infty$ of the sizes of the largest families living at time
$t$, and of the ages of the oldest families living at time $t$. Before this, we give
estimates on the second factorial moment of the number of large families, used repeatedly in
the sequel.

\subsection{A preliminary lemma}
\label{sec:prelim-lem}

Let us recall that, under $\PP_t$, we can adopt the representation of the genealogy at
time $t$ by the coalescent point process $H_0, H_1, H_2,\ldots$, where $H_0=+\infty$ and the
$(H_i;i\ge 1)$ are i.i.d., killed at their first value ($=N_t$) larger than $t$. For all
$i\geq 0$, we call \emph{branch $i$} the lineage represented by $H_i$.

For all $t>0$, $x\geq 1$, $0\leq s_1<s_2\leq +\infty$, we define $K_t(x,s_1,s_2)$ as the
number of haplotypes carried by more than $x$ individuals alive at time $t$, whose original
mutation occurred on the ancestral lineage (branch $0$), and has age in $(s_1, s_2]$ (or,
equivalently, in $(s_1,s_2\wedge t]$).

\begin{lem}
  \label{lem:summary-technique}
  For all $t>0$, $x\geq 1$, $0\leq s_1<s_2\leq t$, we have
  \begin{multline}
    \EE_t [M_t(x,s_1,s_2)(M_t(x,s_1,s_2)-1)] \le 2\ \EE_t
    K_t(x,s_1,s_2)\ \left[1+\frac{\frac{1}{W(s_1)}-\frac{1}{W(t)}}{1-\frac{1}{W(t)}}\ \EE_t
      N_t\right]\times \\
    \times\left[\frac{b(1+\theta(s_2-s_1))}{\alpha}\left(1-\frac{1}{W_\theta(s_2)}\right)^{\lceil
      x\rceil-1}\left(e^{-\theta s_1}+\int_0^{s_1}\theta e^{-\theta
        z}\frac{W_\theta(s_2)-W_\theta(s_1-z)}{W_\theta(s_2)}dz\right)\right. \\
  \left. + 4\ \frac{\frac{1}{W(s_1)}-\frac{1}{W(t)}}{1-\frac{1}{W(t)}}\ \EE_t
  N_t\ \EE_t K_t\ (\lceil x/2\rceil,s_1,s_2)\right]. \label{eq:jolie-equation}
  \end{multline}
  and
  \begin{align}
    \EE_t K_t(x,s_1,s_2) &
    \leq\frac{b}{\alpha}\int_{s_1}^{s_2}\left(1-\frac{1}{W_\theta(y)}\right)^{\lceil
      x\rceil-1}\left(e^{-\theta y}+\int_0^{y}\theta
      e^{-\theta z}\frac{W_\theta(y)-W_\theta(y-z)}{W_\theta(y)}dz\right)(\theta
    dy+\delta_t(dy))  \label{eq:K} \\
    & \leq\frac{b}{\alpha}\int_{s_1}^{s_2}\left(1-\frac{1}{W_\theta(y)}\right)^{\lceil
      x\rceil-1}(\theta
    dy+\delta_t(dy)). \label{eq:K-version-1}
  \end{align}
\end{lem}
This lemma is proved in the Appendix.

\subsection{Convergence in distribution of the size of the most frequent haplotype}

Let us recall the notation $X^{(1)}_t\geq X^{(2)}_t\geq\ldots\geq X^{(k)}_t\geq\ldots$ for the
ordered sequence of sizes of all living families in the population at time $t$ (with the
convention that $X^{(k)}_t=0$ when $k$ is larger than the number of living haplotypes at time
$t$). Our first goal is to prove the convergence in distribution of $X^{(1)}$ using only the
exact formulae~(\ref{eq:exp-FS}) and~(\ref{eq:ancestral-FS}) and the \emph{coalescent point
  process} construction of the genealogy of the splitting tree.

\begin{figure}[t]
  \begin{center}
    \psfrag{0}{0} \psfrag{st}{$s_t$} \psfrag{t}{$t$} \psfrag{F1}{${\cal T}_1$}
    \psfrag{F2}{${\cal T}_2$} \psfrag{F3}{${\cal T}_3$} \psfrag{F4}{${\cal T}_4$}
    \psfrag{F5}{${\cal T}_5$} \psfrag{F6}{${\cal T}_6$} \psfrag{F7}{${\cal T}_7$}
    \psfrag{F8}{${\cal T}_8$} \psfrag{F9}{${\cal T}_9$} \psfrag{F10}{${\cal T}_{10}$}
    \includegraphics[width=15cm]{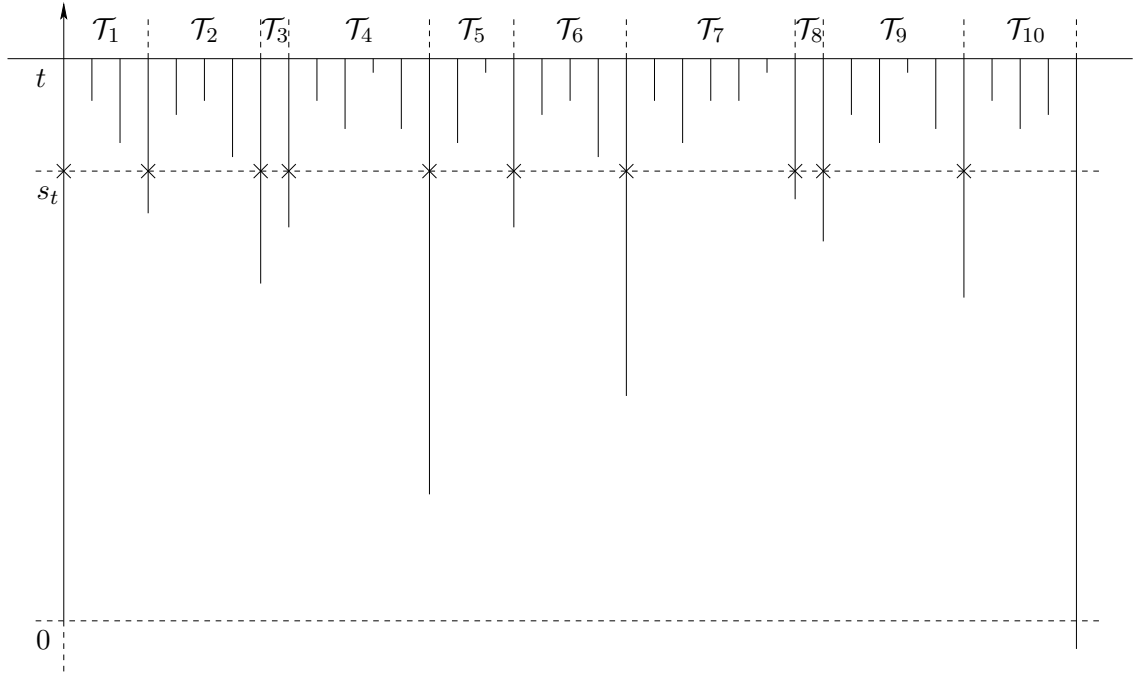}
  \end{center}
  \caption{The definition of the sub-trees $({\cal T}_i)_{1\leq i\leq N_{t,s_t}}$. The first
    vertical line represents the ancestral lineage (branch 0) and the other vertical lines
    have i.i.d.\ lengths $H_1,H_2,\ldots$ The rightmost vertical line is the first one higher
    than $t$, with length $H_{N_t}$. The crosses represent the $N_{t,s_t}$ individuals alive
    at time $s_t$ and having descendants at time $t$. Here, $N_{t,s_t}=10$.}
  \label{fig:cv-loi}
\end{figure}

The general idea is the following one. We divide the population at time $t$ into several
sub-populations corresponding to distinct ancestors at a given time $s$, as shown in
Fig.~\ref{fig:cv-loi}. This gives a sequence of sub-trees $({\cal T}_i)_{1\leq i\leq
  N_{t,s}}$, where $N_{t,s}$ is the number of individuals alive at time $s$ having descendants
at time $t$. These individuals are represented by crosses in Fig.~\ref{fig:cv-loi}. We choose
$s=s_t$ in such a way that $N_{t,s_t}\to\infty$ and the event under consideration (here, the
event that there exists a haplotype carried by more than $x_t$ individuals at time $t$) has a
small probability in each sub-tree and are ``nearly independent'' (in a sense specified in the
proofs below) in distinct sub-trees. The key argument of the proof, for which
Lemma~\ref{lem:summary-technique} is needed, consists in checking that, in each subtree, the
unknown probability that there exists a haplotype satisfying the property under consideration
(here, carried by more than $x_t$ individuals) is close to the expected number of such
haplotypes, which is known explicitly. Here, this reads
$$
\PP_t[L_{t-s_t}(x_t)\geq 1]\sim \EE_t L_{t-s_t}(x_t).
$$

\begin{thm}
  \label{thm:cv-loi-most-frequent}
  Assume $\alpha<\theta$. For all $c\in\RR$, let
  \begin{equation}
    \label{eq:def-x_t-loi}
    x_t(c):=\frac{\alpha t-\frac{\theta}{\theta-\alpha}\log
      t}{|\log\varphi(\theta)|}+c.
  \end{equation}
  Then
  \begin{equation*}
    \PP_t(X^{(1)}_t< x_t(c)) 
    \sim \frac{1}{1+A(\theta)\varphi(\theta)^{c-1+\{-x_t(c)\}}}  
  \end{equation*}
  as $t\rightarrow+\infty$, where $A(\theta)$ is defined in~(\ref{eq:def-A}).
\end{thm}

\paragraph{Proof.}
Set 
$$
F(t,x):=\PP_t(X^{(1)}_t\geq x)=\PP_t[L_t(x)\geq 1],\qquad G(t,x):=\EE_t[L_t(x)]
$$
and
$$
F(t,x,s):=\PP_t[M_t(x,s)\geq 1]. 
$$
We have for all $s<t$
\begin{multline*}
  0\leq F(t,x)-\PP_t\Big(\exists i\in\{1,\ldots,N_{t,s}\}:{\cal T}_i\text{\ contains a haplotype} \\
  \text{\ carried by more than $x$ individuals}\Big)\leq F(t,x,t-s).
\end{multline*}
which also reads
\begin{multline}
  0\leq F(t,x)-1+\EE_t[(1-F(t-s,x))^{N_{t,s}}] \\ =F(t,x)-
  \frac{1}{1+\PP(H>t\mid H>t-s)\,\left(\frac{1}{F(t-s,x)}-1\right)}
  \leq F(t,x,t-s), \label{eq:est-proba}
\end{multline}
since $N_{t,s}$ is a geometric random variable of parameter $\PP(H>t\mid H>t-s)$. In view of
this, in order to find a non-trivial limit for $F(t,x_t)$, we need to find $s_t$ and $x_t$
such that $F(t,x_t,t-s_t)=o(1)$ and $F(t-s_t,x_t)=o(1)$ and is asymptotically equivalent to
\begin{equation}
  \label{eq:equiv-geom-N_t,s}
  \PP(H>t\mid H>t-s_t)=W(t-s_t)/W(t)\sim e^{-\alpha s_t}  
\end{equation}
as $t\rightarrow+\infty$. In order to find an explicit asymptotic equivalent of
$F(t-s_t,x_t)$, we will compare it with $G(t-s_t,x_t)$.

Let us check that the choice $x_t(c)$ in~(\ref{eq:def-x_t-loi}) for $x_t$ and
$$
s_t=s_t(b)=t-b\log t
$$
satisfy the above properties for all $b>1/(\theta-\alpha)$. 

On the one hand, the fact that $F(t,x_t(c),t-s_t(b))=o(1)$ is an immediate consequence of
Proposition~\ref{prop:frequent-a<t}, since $b\log t>(1+\varepsilon)\log t/(\theta-\alpha)$ for
some $\varepsilon>0$. On the other hand, we can compute an asymptotic equivalent of
$G(t-s_t(b),x_t(c))$ following closely the computation of the proof of
Proposition~\ref{prop:frequent-a<t}. Here are the main steps of the computation: we have
\begin{align*}
  G(t-s_t(b),x_t(c)) & =W(t-s_t(b))\left(1-\frac{\psi(\theta)}{\theta}\right)^{\lceil x_t(c)\rceil-1}
  \int_0^{t-s_t(b)}\frac{e^{-\theta x}}{W_\theta(x)}(1-\rho(x))^{\lceil x_t(c)\rceil-1}(\theta
  dx+\delta_{t-s_t(b)}(dx)) \\ &
  \sim\frac{e^{-\alpha s_t(b)}\,t^{\theta/(\theta-\alpha)}\,\varphi(\theta)^{c-1+\{-x_(c)\}}}
  {\psi'(\alpha)} \int_0^{t-s_t(b)}\frac{e^{-\theta x}}{W_\theta(x)}
  (1-\rho(x))^{\lceil x_t(c)\rceil-1}(\theta
  dx+\delta_{t-s_t(b)}(dx)).
\end{align*}
It is easy to deduce from Lemma~\ref{lem:W_theta}~(ii) that the contribution of the Dirac mass
and of the integral on the interval $[0,(1-\varepsilon)\log t/(\theta-\alpha)]$ for any fixed
$\varepsilon\in(0,1)$ are both $o(e^{-\alpha s_t(b)})$. Using the fact that
$W_\theta(x)\rightarrow \theta/\psi(\theta)$ when $x\rightarrow+\infty$ and the change of
variable $y=B(\lceil x_t(c)\rceil-1)e^{-(\theta-\alpha)x}$, the contribution of the integral
on the interval $[(1-\varepsilon)\log t/(\theta-\alpha),b\log t]$ is asymptotically equivalent
to
\begin{multline*}
  \frac{\psi(\theta)e^{-\alpha
      s_t(b)}\varphi(\theta)^{c-1+\{-x_t(c)\}}}{\psi'(\alpha)(\theta-\alpha)}
  \left(\frac{t}{B(\lceil x_t(c)\rceil-1)}\right)^{\theta/(\theta-\alpha)} \\
  \int_{B(\lceil
    x_t(c)\rceil-1)\ t^{-b(\theta-\alpha)}}^{B(\lceil
    x_t(c)\rceil-1)\ t^{-(1-\varepsilon)}} y^{\alpha/(\theta-\alpha)}
  \exp\left((\lceil x_t(c)\rceil-1)\log\left(1-\rho\left(\frac{\log\frac{B(\lceil
            x_t(c)\rceil-1)}{y}}{\theta-\alpha}\right)\right)\right)dy.  
\end{multline*}
As in the proof of Proposition~\ref{prop:frequent-a<t}, Lebesgue's dominated convergence
theorem then yields
\begin{equation}
  \label{eq:youpi} 
  G(t-s_t(b),x_t(c)) \sim A(\theta)\,\varphi(\theta)^{c-1+\{-x_t(c)\}}\PP(H>t\mid H>t-s_t(b))
\end{equation}
as $t\rightarrow+\infty$, where we used the fact that $x_t(c)t^{-b(\theta-\alpha)}\rightarrow
0$ when $t\rightarrow+\infty$ since $b>1/(\theta-\alpha)$.

Now, it only remains to check that 
\begin{equation}
  \label{eq:goal-1}
  G(t-s_t(b),x_t(c))\sim F(t-s_t(b),x_t(c))  
\end{equation}
when $t\rightarrow+\infty$. Since for all $t,x>0$
$$
0\leq G(t,x)-F(t,x)=\EE_t (L_{t}(x))-\PP_t(L_{t}(x)\geq 1)\leq
\EE_t[L_{t}(x)(L_{t}(x)-1)],
$$
it is sufficient to prove that
\begin{equation}
  \label{eq:more-than-2-haplo}
  \EE_t\left[L_{t-s_t(b)}(x_t(c))\left(L_{t-s_t(b)}(x_t(c))-1\right)\right]=o\left(G(t-s_t(b),x_t(c))\right)
  =o(e^{-\alpha s_t(b)}).
\end{equation}
Taking $s_1=0$ and $s_2=t-s_t(b)$ in
Lemma~\ref{lem:summary-technique}~(\ref{eq:jolie-equation}) yields
\begin{multline*}
  \EE_t\left[L_{t-s_t(b)}(x_t(c))\left(L_{t-s_t(b)}(x_t(c))-1\right)\right]\leq
  2\EE_t K_{t-s_t(b)}(x_t(c))(1+\EE_t
  N_{t-s_t(b)})\times \\ \times\left[\frac{b(1+\theta
  (t-s_t(b)))}{\alpha}\left(1-\frac{1}{W_\theta(t-s_t(b))}\right)^{\lceil x_t(c)\rceil-1}+4\EE_t
  N_{t-s_t(b)}\EE_t K_{t-s_t(b)}(x_t(c)/2)\right],
\end{multline*}
where we used the notation $K_s(x):=K_s(x,0,s)$. Using the inequality $\EE_t
K_{s}(x)\leq \EE_t L_s(x)=G(s,x)$ and the estimates
\begin{gather}
  \EE_t N_{t-s_t(b)}=\frac{1}{\PP(H> t-s_t(b))}=W(t-s_t(b))\sim
  \frac{e^{\alpha(t-s_t(b))}}{\psi'(\alpha)}, \label{eq:equiv-N_t,s} \\
  \left(1-\frac{1}{W_\theta(t-s_t(b))}\right)^{\lceil
    x_t(c)\rceil-1}\leq \varphi(\theta)^{\lceil x_t(c)\rceil-1}\leq Ce^{-\alpha
    t}t^{\theta/(\theta-\alpha)}, \notag
\end{gather}
and
\begin{equation*}
  G(t-s_t(b),x_t(c)/2)\sim 2^{\theta/(\theta-\alpha)}\,A(\theta)\,
  \varphi(\theta)^{c/2-1+\{-x_t(c)/2\}}\,
  e^{-\alpha(s_t(b)-t/2)}\,t^{-\theta/2(\theta-\alpha)},
\end{equation*}
which can be proved following the same computation as above, we finally obtain for $t$ large
enough
\begin{multline*}
  \EE_t\left[L_{t-s_t(b)}(x_t(c))\left(L_{t-s_t(b)}(x_t(c))-1\right)\right] \\ \leq
  C\,G(t-s_t(b),x_t(c))\left(e^{\alpha(t-s_t(b))}\,\log t\,e^{-\alpha t}t^{\theta/(\theta-\alpha)}+
  e^{2\alpha(t-s_t(b))}\,e^{-\alpha(s_t(b)-t/2)}\,t^{-\theta/2(\theta-\alpha)}\right),
\end{multline*}
which entails~(\ref{eq:more-than-2-haplo}) and concludes the proof of
Theorem~\ref{thm:cv-loi-most-frequent}.\hfill$\Box$
\bigskip

\subsection{Joint convergence in law of the sizes of the most abundant families}

The general idea of the previous argument can be summarized as follows: we construct a random number
$N_{t,s_t}$ of i.i.d.\ random variables $X_1,\ldots,X_{N_t,s_t}$---the sizes of the most frequent haplotype in
each sub-tree---such that $X^{(1)}_t=\max\{X_1,\ldots,X_{N_{t,s_t}}\}$ with high probability. Our previous
result then corresponds to a classical argument of extreme value theory, which is known to extend easily to
compute the extremes statistics and the joint law of the largest random variables among $X_1,\ldots,X_{N_{t,s_t}}$. The object of this
subsection is to prove that this argument is valid in our situation.

\begin{thm}
  \label{thm:freq-2-t<a}
  Assume $\alpha<\theta$ and recall the definition of $x_t(c)$ in~(\ref{eq:def-x_t-loi}). For all $n\in\NN$,
  $k_1,\ldots,k_n\in\ZZ_+$ and $c_1,\ldots,c_n$ such that $c_i\geq c_{i+1}+1$ for all $i\in\{1,\ldots,n-1\}$,
  we have, as $t\to +\infty$,
  \begin{multline}
  \PP_t\Big[L_t(x_t(c_1))=k_1,\ L_t(x_t(c_2))-L_t(x_t(c_1))=k_2,\ldots,\
    L_t(x_t(c_n))-L_t(x_t(c_{n-1}))=k_n\Big] \\
  \sim \binom{k_1+\ldots+k_n}{k_1,\ldots,k_n}\,
  \frac{\tau_t(c_1)^{k_1}(\tau_t(c_2)-\tau_t(c_1))^{k_2}\ldots(\tau_t(c_n)-\tau_t(c_{n-1}))^{k_n}}
  {(1+\tau_t(c_n))^{k_1+\ldots+k_n+1}},   \label{eq:goal-2}
  \end{multline}
  with
  $$
  \tau_t(c):= A(\theta)\varphi(\theta)^{c-1+\{-x_t(c)\}}
  $$
  for all $c\in\RR$, where the constant $A(\theta)$ is defined in~(\ref{eq:def-A}).
\end{thm}

\paragraph{Proof.}
Let us denote by $A(t;c_1,\ldots,c_n;k_1,\ldots,k_n)$ the event in the probability in the l.h.s.\
of~(\ref{eq:goal-2}). Using the notation of the proof of Theorem~\ref{thm:cv-loi-most-frequent}, for all
$b>1/(\theta-\alpha)$, we define $B(t,b;c_1,\ldots,c_n;k_1,\ldots,k_n)$ the event that among the sub-trees
${\cal T}_1,\ldots,{\cal T}_{N_t,s_t(b)}$, there are exactly $k_i$ haplotypes carried by a number of living
individuals at time $t$ belonging to $[x_t(c_{i+1}),x_t(c_i))$ for all $0\leq i\leq n-1$, with the convention
$c_0=+\infty$. Then
\begin{multline*}
  \Big|\PP_t(A(t;c_1,\ldots,c_n;k_1,\ldots,k_n))-\PP_t(B(t,b;c_1,\ldots,c_n;k_1,\ldots,k_n)\Big| \\
  \leq\PP_t[M_t(x_t(c_n),t-s_t(b))\geq 1]= F(t,x_t(c_n),t-s_t(b))=o(1).
\end{multline*}
Now, for all fixed $t\geq 0$, using the notation $X_1,\ldots,X_{N_t,s_t}$ introduced above, we define for all
$a\in\RR$
$$
S_t(a)=\#\{i\leq N_{t,s_t}:X_i\geq a\}.
$$
Defining
\begin{multline*}
  C(t,b;c_1,\ldots,c_n;k_1,\ldots,k_n):= \\ \Big\{S_t(x_t(c_1))=k_1,\ S_t(x_t(c_2))-S_t(x_t(c_1))=k_2,\ldots,\
  S_t(x_t(c_n))-S_t(x_t(c_{n-1}))=k_n\Big\},
\end{multline*}
we have
\begin{multline*}
  \Big|\PP_t(B(t,b;c_1,\ldots,c_n;k_1,\ldots,k_n))-\PP_t(C(t,b;c_1,\ldots,c_n;k_1,\ldots,k_n))\Big| \\
  \begin{aligned}
    & \leq\PP_t\Big[\exists i\leq N_{t,s_t(b)}:{\cal T}_i\text{\ contains at least 2 haplotypes carried by
      more than $x_t(c_n)$ individuals}\Big] \\
    & \leq\sum_{k\geq 1}\PP(H> t\mid H>t-s_t(b))\,\PP(H\leq t\mid H>
    t-s_t(b))^{k-1}\, k\, \EE_t\Big[L_{t-s_t(b)}(x_t(c_n))(L_{t-s_t(b)}(x_t(c_n))-1)\Big] \\
    & =\frac{\EE_t\left[L_{t-s_t(b)}(x_t(c_n))(L_{t-s_t(b)}(x_t(c_n))-1)\right]}{\PP(H>t\mid
      H>t-s_t(b))},
  \end{aligned}
\end{multline*}
where we used the fact that $N_{t,s_t(b)}$ has geometric distribution with parameter $\PP(H> t\mid
H>t-s_t(b))$. Equations~(\ref{eq:equiv-geom-N_t,s}) and~(\ref{eq:more-than-2-haplo}) then yield
\begin{equation*}
  \Big|\PP_t(A(t;c_1,\ldots,c_n;k_1,\ldots,k_n))-\PP_t(C(t,b;c_1,\ldots,c_n;k_1,\ldots,k_n)\Big|=o(1).
\end{equation*}

Next, $\PP_t(C(t,b;c_1,\ldots,c_n;k_1,\ldots,k_n))$ can be computed using standard extreme value
techniques: conditioning on $N_{t,s_t(b)}$ and considering all the possible ways to realize this event, we
have
\begin{multline*}
  \PP_t(C(t,b;c_1,\ldots,c_n;k_1,\ldots,k_n)) \\ 
  \begin{aligned}
    & =\sum_{k\geq k_1+\ldots+k_n}\PP(H> t\mid H>t-s_t(b))\,\PP(H\leq t\mid H>
    t-s_t(b))^{k-1}\,\binom{k}{k_1,\ldots,k_n}\times \\ & \qquad \PP_t\Big[X_1,\ldots,X_{k_1}\geq
    x_t(c_1)>X_{k_1+1},\ldots, X_{k_1+k_2}\geq x_t(c_2)>\ldots\geq
    x_t(c_{n})>X_{k_1+\ldots+k_{n}+1},\ldots,X_k\Big] \\ & =\sum_{k\geq k_1+\ldots+k_n}
    \PP(H> t\mid H>t-s_t(b))\,\PP(H\leq t\mid H>
    t-s_t(b))^{k-1}\,\binom{k_1+\ldots+k_n}{k_1,\ldots,k_n}\,\binom{k}{k_1+\ldots+k_n}\times \\ & \qquad\qquad
    F(t-s_t(b),x_t(c_1))^{k_1}\,
    [F(t-s_t(b),x_t(c_2))-F(t-s_t(b),x_t(c_1))]^{k_2}\times\ldots \\ & \qquad\qquad\qquad
    \ldots\times[F(t-s_t(b),x_t(c_n))-F(t-s_t(b),x_t(c_{n-1}))]^{k_n}\,[1-F(t-s_t(b),x_t(c_n))]^{k-k_1-\ldots-k_n}.
  \end{aligned}
\end{multline*}
The equation
$$
\sum_{k\geq m}\binom{k}{m}x^{k-m}=\frac{1}{(1-x)^{m+1}},\qquad\forall x\in(-1,1),\ m\in\NN
$$
yields
\begin{multline*}
  \PP_t(C(t,b;c_1,\ldots,c_n;k_1,\ldots,k_n)) \\ 
  =\frac{\PP(H>t\mid H>t-s_t(b))\PP(H\leq t\mid H>t-s_t(b))^{k_1+\ldots+k_n-1}}
  {\left[\PP(H>t\mid H>t-s_t(b))+F(t-s_t(b),x_t(c_n))\PP(H\leq t\mid H>t-s_t(b))\right]^{k_1+\ldots+k_n+1}}\times \\
  \binom{k_1+\ldots+k_n}{k_1,\ldots,k_n}\,F(t-s_t(b),x_t(c_1))^{k_1}\,
    [F(t-s_t(b),x_t(c_2))-F(t-s_t(b),x_t(c_1))]^{k_2}\times\ldots \\
    \ldots\times[F(t-s_t(b),x_t(c_n))-F(t-s_t(b),x_t(c_{n-1}))]^{k_n},
\end{multline*}
and Theorem~\ref{thm:freq-2-t<a} then follows from~(\ref{eq:equiv-geom-N_t,s}),~(\ref{eq:youpi})
and~(\ref{eq:goal-1}).\hfill$\Box$ 
\bigskip

In order to state our next result, we define the number $t_n$ by the equation $x_{t_n}(0)=n$. In view
of~(\ref{eq:def-x_t-loi}), this equation has a unique solution $t_n$ if $n$ is large enough, say larger than
$n_0$. In addition,
$$
t_n\sim\frac{|\log\varphi(\theta)|\,(\theta-\alpha)}{\theta}\,n,
$$
and hence $t_n\rightarrow+\infty$ as $n\rightarrow+\infty$.

We also recall the notation ${\cal M}(\RR)$ for the set of nonnegative $\sigma$-finite measures on $\RR$, finite on $\RR_+$,
and the definition of the \emph{semi-vague} topology as the one induced by all maps of the form
$$
\nu\in{\cal M}(\RR)\mapsto \int_{\RR}u(x)\nu(dx),
$$
for all continuous bounded function $u$ on $\RR$ such that there exists $x_0\in\RR$ such that $u(x)=0$ for all
$x\leq x_0$. Note that this topology is stronger than the usual vague topology, but weaker than the usual weak
topology.

\begin{cor}
  \label{cor:freq-a<t}
  Assume $\alpha<\theta$. Then, the sequence of point processes $(Z_n)_{n\geq n_0}$ on $\ZZ$, defined by
  $$
  Z_n:=\sum_{k\geq 1}\delta_{X^{(k)}_{t_n}-n},
  $$
  converges as $n\rightarrow+\infty$ in $\PP^\star$-distribution on set ${\cal M}(\RR)$ equipped with the semi-vague topology to a
  mixed Poisson point measure on $\ZZ$ with intensity measure
  $$
  {\cal E}
  A(\theta)\frac{\psi(\theta)}{\theta}\,\sum_{c\in\ZZ}\varphi(\theta)^{c-1}\,\delta_c,
  $$
  where the mixture coefficient $\cal E$ has exponential distribution with parameter $1$.
\end{cor}


The proof of such results is quite standard (cf.~\cite{Kallenberg} in the general context of
point processes and~\cite{Leadbetter} more specifically on extreme values). However, we shall
give a proof for sake of completeness and because of the specificity of the semi-vague
topology.

Note that one can also easily obtain the convergence, in the sense of finite-dimensional
distributions, of any finite sequence
of translated extreme family sizes towards the corresponding sequence of extreme points of the limit
point process in the previous result. We shall not prove this, but instead we refer
to~\cite{Leadbetter} for the proof of similar standard results.

\paragraph{Proof.}
Let us first prove the convergence in distribution when ${\cal M}(\RR)$ is equipped with the vague
topology. This amounts to prove the joint convergence in distribution of the random variables
$L_{t_n}(n+i)-L_{t_n}(n+i+1))=L_{t_n}(x_{t_n}(i))-L_{t_n}(x_{t_n}(i+1))$, giving the number of haplotypes represented by exactly $n+i$ individuals at time $t_n$, for $b\leq i\leq a$ for all $a<b$ in $\ZZ$. Fix $b<a$ in
$\ZZ$ and $k_b,\ldots,k_a$ in $\ZZ_+$. On the one hand, we claim that
\begin{multline*}
  \lim_{n\rightarrow+\infty}\PP_{t_n}(L_{t_n}(n+i)-L_{t_n}(n+i+1)=k_i,\ \forall b\leq i\leq a) \\
  \begin{aligned}
    & =\lim_{n\rightarrow+\infty}\sum_{k\geq 0}\PP_{t_n}(L_{t_n}(x_{t_n}(i))-L_{t_n}(x_{t_n}(i+1))=k_i,\
    \forall b\leq i\leq a\text{\ and\ }L_{t_n}(x_{t_n}(a+1))=k) \\
    & =\sum_{k\geq 0}\binom{k_a+\ldots+k_b+k}{k_a,\ldots,k_b,k}\,
    \frac{\tau_{t_n}(a+1)^{k}(\tau_{t_n}(a)-\tau_{t_n}(a+1))^{k_a}\ldots(\tau_{t_n}(b)-\tau_{t_n}(b+1))^{k_b}}
    {(1+\tau_t(b))^{k+k_b+\ldots+k_a+1}} \\ &
    =\binom{k_a+\ldots+k_b}{k_a,\ldots,k_b}\,\frac{[A(1-\varphi)]^{k_a+\ldots+k_b}\,\varphi^{(a-1)k_a+\ldots+(b-1)k_b}}
    {[1+A(\varphi^{b-1}-\varphi^a)]^{k_a+\ldots+k_b+1}},
  \end{aligned}
\end{multline*}
where we used the fact that, for all $x\in\ZZ$ and $n\geq n_0$,
$$
\tau_{t_n}(x)=A \varphi^{x-1}\qquad \text{with}\qquad A=A(\theta) \quad \mbox{and}\quad\varphi=\varphi(\theta):=1-\frac{\psi(\theta)}{\theta}.
$$
This is an immediate consequence of Theorem~\ref{thm:freq-2-t<a}, provided we can justify the exchange of the
sum over $k$ and the limit $n\rightarrow+\infty$, i.e.\ that we can control the remainder of the series
uniformly over $n\geq n_0$. The following inequality, making use of Proposition~\ref{prop:frequent-a<t},
solves this question: for all $N\in\NN$,
\begin{multline*}
  \PP_{t_n}(L_{t_n}(x_{t_n}(i))-L_{t_n}(x_{t_n}(i+1))=k_i,\ \forall b\leq i\leq a\text{\ and\
  }L_{t_n}(x_{t_n}(a+1))\geq N) \\ \leq\PP_{t_n}(L_{t_n}(x_{t_n}(a+1))\geq N)\leq\frac{\EE_{t_n}
    L_{t_n}(x_{t_n}(a+1))}{N}\leq\frac{1}{N}\sup_{t\geq 0}\EE_{t} L_t(x_t(a+1))\leq\frac{C}{N}
\end{multline*}
for some constant $C>0$.

On the other hand, assume that $\cal E$ is an exponential random variable with parameter 1,
and $(P_x)_{x\in\ZZ}$ is a sequence of r.v.\ with $P_x$ distributed as a mixed Poisson with
parameter ${\cal E} A(1-\varphi)\varphi^{x-1}$ and such that the r.v.\ $(P_x)_{x\in\ZZ}$ are
independent conditionally on ${\cal E}$. Then,
\begin{align*}
  \PP(P_i=k_i,\ \forall b\leq i\leq a)
  & =\int_0^{\infty}dx\,e^{-x}\,e^{-A(1-\varphi)x\sum_{m=b}^a\varphi^{m-1}}
  \frac{(xA(1-\varphi))^{k_a+\ldots+k_b}\,\varphi^{(a-1)k_a+\ldots+(b-1)k_b}}{k_a!\ldots k_b!} \\
  & =\frac{[A(1-\varphi)]^{k_a+\ldots+k_b}\,\varphi^{(a-1)k_a+\ldots+(b-1)k_b}}{k_a!\ldots k_b!\,
    [1+A(1-\varphi)\sum_{m=b}^a\varphi^{m-1}]^{k_a+\ldots+k_b+1}}\,\int_0^{\infty}y^{k_a+\ldots+k_b}e^{-y}dy \\
  & =\binom{k_a+\ldots+k_b}{k_a,\ldots,k_b}\,\frac{[A(1-\varphi)]^{k_a+\ldots+k_b}\,\varphi^{(a-1)k_a+\ldots+(b-1)k_b}}
  {[1+A(\varphi^{b-1}-\varphi^a)]^{k_a+\ldots+k_b+1}},
\end{align*}
where we used the change of variable $y=x(1+A(1-\varphi)\sum_{m=b}^a\varphi^{m-1})$. Observing
that $\PP_{t_n}$ converges to $\PP^\star$ for the total variation norm, this ends the proof of
Corollary~\ref{cor:freq-a<t} when ${\cal M}(\RR)$ is equipped with the vague topology.

To complete the proof of Corollary~\ref{cor:freq-a<t}, since all the point measures $Z_n$ have
support in $\ZZ$, we need to check that, for any continuous bounded function $f$ on $\RR$ and
any sequence $(u(k))_{k\in\ZZ}$ such that $u(k)=0$ for all $k\leq k_0$ for some $k_0\in\ZZ$,
$$
\lim_{n\rightarrow+\infty}\EE_{t_n}\,f\left(\int_{\ZZ}u(x) Z_n(dx)\right)=\EE\, f\Bigg(\sum_{k>k_0}u(k)P_k\Bigg).
$$
Note first that the sum in the r.h.s.\ is almost surely finite since, conditional on ${\cal E}$, this is a sum of
independent r.v.\ with only finitely many of them being non-zero by Borel-Cantelli's lemma.

Next, fix $\varepsilon>0$ and let $a$ and $T$ be large enough so that
\begin{equation*}
  \frac{\tau_{t_n}(a)}{1+\tau_{t_n}(a)}=\frac{A\varphi^{a-1}}{1+A\varphi^{a-1}}\leq A\varphi^{a-1}\leq\varepsilon,
\end{equation*}
and, by Theorem~\ref{thm:cv-loi-most-frequent},
$$
\sup_{t\geq T}\,\PP_t \left(X^{(1)}_t\geq x_t(a)\right)\leq 2\varepsilon.
$$
Then, for all $n$ such that $t_n\geq T$,
$$
\Bigg|\EE_{t_n}\, f\left(\int_\ZZ u(x)Z_n(dx)\right)
-\EE_{t_n}\, f\Bigg(\sum_{k=k_0+1}^a u(k)\big(L_{t_n}(x_{t_n}(k))-L_{t_n}(x_{t_n}(k+1))\big)\Bigg)\Bigg|
\leq 4\varepsilon\|f\|_\infty.
$$
The first step of the proof then yields
$$
\Bigg|\EE_{t_n}\, f\left(\int_\ZZ u(x)Z_n(dx)\right)
-\EE\, f\Bigg(\sum_{k=k_0+1}^a u(k)P_k\big)\Bigg)\Bigg|
\leq (4\|f\|_\infty+1)\varepsilon
$$
for all $n$ large enough. Since
$$
\PP(\exists k\geq a:P_k\geq
1)\leq\sum_{k\geq a}\PP(P_k\geq 1)=\sum_{k\geq a}(1-e^{-A\varphi^{k-1}})\leq A\sum_{k\geq a}\varphi^{k-1}\leq
\frac{\varepsilon}{1-\varphi},
$$
we finally obtain
$$
\Bigg|\EE_{t_n}\, f\left(\int_\ZZ u(x)Z_n(dx)\right)
-\EE\, f\Bigg(\sum_{k>k_0} u(k)P_k\big)\Bigg)\Bigg|
\leq (4\|f\|_\infty+1+(1-\varphi)^{-1})\varepsilon,
$$
which ends the proof of Corollary~\ref{cor:freq-a<t}.\hfill$\Box$
\bigskip

\subsection{Convergence in distribution of the ages of oldest families}

The previous method can  easily be extended to prove the convergence of the ages of the oldest families. Let
us recall the notation $A^{(1)}_t\geq A^{(2)}_t\geq\ldots\geq A^{(k)}_t\geq \ldots$ for the ordered
sequence of ages of all alive families at time $t$ (with the convention that
$A^{(k)}_t=0$ when $k$ is larger than the number of alive families at time $t$).


\begin{thm}
  \label{thm:cv-distr-old-a<t}
  Assume $\alpha<\theta$ and define for all $a\in\RR$
  $$
  x_t(a)=\frac{\alpha t}{\theta}+a.
  $$
  For all $n\in\NN$, $k_1,\ldots,k_n\in\ZZ_+$ and $a_1>a_2>\ldots>a_n$, we have
  \begin{multline}
    \lim_{t\rightarrow+\infty}\,\PP_t\Big[O_t(x_t(a_1))=k_1,\ O_t(x_t(a_2))-O_t(x_t(a_1))=k_2,\ldots,\
    O_t(x_t(a_n))-O_t(x_t(a_{n-1}))=k_n\Big] \\
    =\frac{\theta\psi'(\alpha)}{\psi(\theta)}\,\binom{k_1+\ldots+k_n}{k_1,\ldots,k_n}\,
    \frac{e^{-\theta k_1 a_1}(e^{-\theta a_2}-e^{-\theta a_1})^{k_2}\ldots(e^{-\theta a_n}-e^{-\theta a_{n-1}})^{k_n}}
    {\left(\frac{\theta\psi'(\alpha)}{\psi(\theta)}+e^{-\theta a_n}\right)^{k_1+\ldots+k_n+1}}. \label{eq:goal-old}
  \end{multline}
  In addition, the family of ${\cal M}(\RR)$-valued random variables $(Z_t,t\geq 0)$, defined for all $t\geq
  0$ by
  $$
  Z_t:=\sum_{k\geq 1}\delta_{A^{(k)}_t-\frac{\alpha t}{\theta}},
  $$
  converges as $t\rightarrow+\infty$ in $\PP^\star$-distribution in ${\cal M}(\RR)$ equipped with the semi-vague topology to a mixed
  Poisson point measure on $\RR$ with intensity measure
  \begin{equation}
    \label{eq:intens-measure-P-age}
    {\cal E}\frac{\psi(\theta)}{\psi'(\alpha)}\,e^{-\theta a}\,da,    
  \end{equation}
  where the mixture coefficient ${\cal E}$ has exponential distribution with parameter $1$.
\end{thm}

The proof follows closely the lines of those of Theorems~\ref{thm:cv-loi-most-frequent}
and~\ref{thm:freq-2-t<a} and Corollary~\ref{cor:freq-a<t}. As a first step, we prove the following lemma.

\begin{lem}
  \label{lem:age}
  With the same notation as in Theorem~\ref{thm:cv-distr-old-a<t}, we have for all $a\in\RR$
  $$
  \lim_{t\rightarrow+\infty}\,\PP_t\left(A^{(1)}_t\leq \frac{\alpha
      t}{\theta}+a\right)=\frac{1}{1+\frac{\psi(\theta)}{\theta\psi'(\alpha)}\,e^{-\theta a}}.
  $$
\end{lem}

\paragraph{Proof.}

The proof of this lemma is similar to the one of Theorem~\ref{thm:cv-loi-most-frequent}. Defining for all
$t,x\geq 0$
$$
F(t,x)=\PP_t(A^{(1)}_t\geq x)=\PP_t[O_t(x)\geq 1]\qquad\text{and}\qquad G(t,x)=\EE_t[O_t(x)],
$$
we have for all $x<s<t$
\begin{align}
  0 & \leq F(t,x)-\PP_t\Big(\exists i\in\{1,\ldots,N_{t,s}\}:{\cal T}_i\text{\ contains at time $t$ a
    haplotype older than $x$}\Big) \notag \\ & = F(t,x)-\frac{1}{1+\PP(H>t\mid H>t-s)\left(\frac{1}{F(t-s,x)}-1\right)} \leq F(t,t-s)\leq G(t,t-s). \label{eq:chouette}
\end{align}
Defining 
$$
s_t(b)=bt\qquad\text{with}\qquad b\in(0,1-\alpha/\theta),
$$
following the proof of Proposition~\ref{prop:age-subcrit}, one easily checks that $G(t,t-s_t(b))=o(1)$ and
\begin{equation}
  \label{eq:cv-distri-age}
  G(t-s_t(b),x_t(a))\sim \frac{\psi(\theta)}{\theta\psi'(\alpha)}\,e^{-\theta a-\alpha s_t(b)},
\end{equation}
where we stick to the notation $x_t(a)=\frac{\alpha t}{\theta}+a$.
Then Lemma~\ref{lem:age} follows from~(\ref{eq:equiv-geom-N_t,s}) and the fact that
$$
G(t-s_t(b),x_t(a))\sim F(t-s_t(b),x_t(a))
$$
as $t\rightarrow\infty$. To prove this last equation, it is sufficient to prove that
\begin{equation}
  \label{eq:goal-calcul-age-loi}
  \EE_t[O_{t-s_t(b)}(x_t(a))(O_{t-s_t(b)}(x_t(a))-1)]=o(e^{-\alpha s_t(b)}), 
\end{equation}
as in \eqref{eq:more-than-2-haplo}.

Now, we observe that
\begin{equation}
  \label{eq:terme-qui-sauve}
  \PP(H>s\mid H\leq
  t)=\frac{\frac{1}{W(s)}-\frac{1}{W(t)}}{1-\frac{1}{W(t)}}\sim \frac{1}{W(s)}\sim  
  \psi'(\alpha)e^{-\alpha s}\quad\text{when\ }s,t\rightarrow+\infty\text{\ with\ }t-s\rightarrow+\infty
\end{equation}
and, by Lemma~\ref{lem:W_theta}~(ii), for a constant $C$ that may change from line to line,
\begin{multline*}
  e^{-\theta x_t(a)}+\int_0^{x_t(a)}\theta e^{-\theta
    z}\frac{W_\theta(t-s_t(b))-W_\theta(x_t(a)-z)}{W_\theta(t-s_t(b))}dz \\
  \begin{aligned}
    & =e^{-\theta x_t(a)}+\int_0^{x_t(a)}\theta e^{-\theta
      z}W_\theta(x_t(a)-z)\varphi(\theta)[\rho(x_t(a)-z)-\rho( t-s_t(b))]dz \\ &
    \leq e^{-\theta x_t(a)}+ C\int_0^{x_t(a)} e^{-\theta
        z}\,e^{-(\theta-\alpha)\left(x_t(a)-z\right)}dz \\ & \leq C e^{-(\theta-\alpha)x_t(a)}.
  \end{aligned}
\end{multline*}
Combining these two facts with~(\ref{eq:equiv-N_t,s}) and
Lemma~\ref{lem:summary-technique}~(\ref{eq:jolie-equation}) in which we take $x=1$,
$s_1=x_t(a)$ and $s_2=t-s_t(b)$, we have
\begin{multline}
  \EE_t\left[O_{t-s_t(b)}(x_t(a))\left(O_{t-s_t(b)}(x_t(a))-1\right)\right] \leq C
  \EE_t K_{t-s_t(b)}(1,x_t(a),t-s_t(b))\left(1+e^{-\alpha
      x_t(a)}\,e^{\alpha(t-s_t(b))}\right)\times \\
  \times\left[(1+t)e^{-(\theta-\alpha)x_t(a)}+e^{-\alpha x_t(a)}\,e^{\alpha(t-s_t(b))}\,
    \EE_t K_{t-s_t(b)}(1,x_t(a),t-s_t(b))\right]. \label{eq:calcul-age-loi}
\end{multline}
Here, in contrast with the proof of Theorem~\ref{thm:cv-loi-most-frequent}, the bound $K_s(1,x,s)\leq
O_s(x)$ is not sufficient to obtain the desired result. Instead, we use Lemma~\ref{lem:summary-technique}~(\ref{eq:K}):
\begin{multline*}
  \EE_t K_{t-s_t(b)}(1,x_t(a),t-s_t(b)) \\ \leq\frac{b}{\alpha}\int_{x_t(a)}^{t-s_t(b)}
  \left(e^{-\theta y}+\int_0^{y}\theta e^{-\theta z} W_\theta(y-z)\varphi(\theta)
  [\rho(y-z)-\rho(y)]dz\right)(\theta dy+\delta_{t-s_t(b)}(dy)).
\end{multline*}
Hence, by Lemma~\ref{lem:W_theta}~(ii) again,
\begin{align*}
  \EE_t K_{t-s_t(b)}(1,x_t(a),t-s_t(b)) & \leq C
  \int_{x_t(a)}^{t-s_t(b)}\left(e^{-\theta y}+\int_0^{y}
    e^{-(\theta-\alpha)(y-z)}e^{-\theta z}dz\right)(\theta dy+\delta_{t-s_t(b)}(dy)) \\ & \leq
  C\int_{x_t(a)}^{t-s_t(b)}e^{-(\theta-\alpha)y}(\theta dy+\delta_{t-s_t(b)}(dy)) \\ & \leq C
  e^{-(\theta-\alpha)x_t(a)}.
\end{align*}
Together with~(\ref{eq:calcul-age-loi}), this yields
$$
\EE_t\left[O_{t-s_t(b)}(x_t(a))\left(O_{t-s_t(b)}(x_t(a))-1\right)\right]\leq C(a)e^{-\alpha s_t(b)}\left(t
  e^{-(1-\alpha/\theta)\alpha t}+e^{-\alpha bt}\right)=o(e^{-\alpha s_t(b)}),
$$
where the constant $C(a)$ depends on $a$ but not on $t$.
This concludes the proof of Lemma~\ref{lem:age}.\hfill$\Box$
\bigskip

\paragraph{Proof of Theorem~\ref{thm:cv-distr-old-a<t}}
Equation~(\ref{eq:goal-old}) can be deduced from Lemma~\ref{lem:age} exactly as Theorem~\ref{thm:freq-2-t<a}
was deduced from Theorem~\ref{thm:cv-loi-most-frequent}. We leave the details to the reader.

In view of~(\ref{eq:goal-old}), the computation in the proof of Corollary~\ref{cor:freq-a<t} immediately
proves (replacing $\tau_{t_n}(a)$ with $\frac{\psi(\theta)}{\theta\psi'(\alpha)}\,e^{-\alpha a}$) that, for
all $a_1>a_2>\ldots>a_n$, the random vector
$$
\Big(O_t(x_t(a_2))-O_t(x_t(a_1)),\ldots,O_t(x_t(a_n))-O_t(x_t(a_{n-1}))\Big)
$$
converges in distribution as $t\rightarrow+\infty$ to a vector whose coordinates are independent
conditionally on $\cal E$ and have mixed Poisson distributions with mixture coefficient $\cal E$ and parameters
$$
{\cal E}\frac{\psi(\theta)}{\theta\psi'(\alpha)}(e^{-\theta a_2}-e^{-\theta
  a_1}),\ldots,{\cal E}\frac{\psi(\theta)}{\theta\psi'(\alpha)}(e^{-\theta a_n}-e^{-\theta a_{n-1}}).
$$

It is then standard to deduce the convergence in distribution of $Z_t$ to $P$ on ${\cal
  M}(\RR)$ equipped with the vague topology (cf.\ e.g.~\cite[Thm.\:4.7]{Kallenberg}). The
semi-vague topology can then be handled similarly as in the proof of
Corollary~\ref{cor:freq-a<t}. Again, we leave the details to the reader.
\hfill$\Box$
\bigskip

\section{Large or old families: convergence in distribution for critical clonal families}
\label{sec:crit}

The method that we used in the previous section can also be applied to the case where $\alpha=\theta$. All the
proofs are similar, and we will only give details at places where the proofs differ. We use the same notation as in the previous section.

\subsection{Frequent haplotypes}
\label{sec:freq-crit}

\begin{thm}
  \label{thm:cv-loi-freq-crit}
  Assume $\alpha=\theta$. For all $c\in\RR$, let
  $$
  x_t(c)=\frac{\alpha^2}{4\psi'(\alpha)}\left(t-\frac{\log t}{2\alpha}+c\right)^2.
  $$
  For all $n\in\NN$, $k_1,\ldots,k_n\in\ZZ_+$ and $c_1>c_2>\ldots>c_n$, we have
  \begin{multline*}
    \lim_{t\rightarrow+\infty}\,\PP_t\Big[L_t(x_t(c_1))=k_1,\ L_t(x_t(c_2))-L_t(x_t(c_1))=k_2,\ldots,\
    L_t(x_t(c_n))-L_t(x_t(c_{n-1}))=k_n\Big] \\
    =\sqrt{\frac{\alpha}{2\pi}}\,e^{-B+\frac{\psi'(\alpha)}{2}}\,\binom{k_1+\ldots+k_n}{k_1,\ldots,k_n}\,
    \frac{e^{-\alpha k_1 c_1}(e^{-\alpha c_2}-e^{-\alpha c_1})^{k_2}\ldots(e^{-\alpha c_n}-e^{-\alpha c_{n-1}})^{k_n}}
    {\left(\sqrt{\frac{\alpha}{2\pi}}\,e^{-B+\frac{\psi'(\alpha)}{2}}+e^{-\alpha
          c_n}\right)^{k_1+\ldots+k_n+1}},
  \end{multline*}
  where the constant $B$ is defined in~(\ref{eq:def-B-crit}). In addition, the family of ${\cal
    M}(\RR)$-valued random variables $(Z_t,t\geq 0)$, defined for all $t\geq 0$ by
  $$
  Z_t:=\sum_{k\geq 1}\delta_{\sqrt{X^{(k)}_t}-\frac{\alpha}{2\sqrt{\psi'(\alpha)}}\left(t-\frac{\log
        t}{2\alpha}\right)},
  $$
  converges as $t\rightarrow+\infty$ in $\PP^\star$-distribution in ${\cal M}(\RR)$ equipped with the semi-vague topology to a mixed
  Poisson point measure on $\RR$ with intensity measure
  \begin{equation}
    \label{eq:intens-measure-P}
    {\cal E}\,\sqrt{2\pi\alpha}\,e^{B-\frac{\psi'(\alpha)}{2}}\,e^{-\alpha c}\,dc,
  \end{equation}
  where the mixture coefficient ${\cal E}$ has exponential distribution with parameter $1$.
\end{thm}
The proof of this result is exactly the same as for Theorems~\ref{thm:freq-2-t<a}
and~\ref{thm:cv-distr-old-a<t}, provided we can prove the following lemma.

\begin{lem}
  \label{lem:crit-freq}
  With the same notation as in Theorem~\ref{thm:cv-loi-freq-crit},  for all $c\in\RR$,
  $$
  \lim_{t\rightarrow+\infty}\
  \PP_t(X_t^{(1)}<x_t(c))=\frac{1}{1+\sqrt{\frac{2\pi}{\alpha}}\,e^{B-\frac{\psi'(\alpha)}{2}}\,e^{-\alpha
      c}}.
  $$
\end{lem}

\paragraph{Proof.}

The proof of this result is similar to the one of Theorem~\ref{thm:cv-loi-most-frequent}. Fix
$\varepsilon>0$. We first observe that Proposition~\ref{prop:a=t} implies that
$$
\PP_t \left(M_t\Big(x_t(c),0,\frac{1-\varepsilon}{2}t\Big)\geq 1\right)\leq
\EE_t M_t\Big(x_t(c),0,\frac{1-\varepsilon}{2}t\Big)=o(1),
$$
and thus
$$
\PP_t(X_t^{(1)}\geq x_t(c))=\PP_t\left(M_t\Big(x_t(c),\frac{1-\varepsilon}{2}t\Big)\geq 1\right)+o(1),
$$
so that it is enough for us to study $\PP_t (M_t(x_t(c),s^{(1)}_t)\geq 1)$, where we put 
$$
s^{(1)}_t:=\frac{1-\varepsilon}{2}t.
$$
Defining
$$
F(t,x)=\PP_t[M_t(x,s^{(1)}_t)\geq 1],\qquad
G(t,x)=\EE_t [M_t(x,s^{(1)}_t)]
$$
and
$$
F(t,x,s)=\PP_t[M_t(x,s)\geq 1], 
$$
we can make the same computation as in the proof of Theorem~\ref{thm:cv-loi-most-frequent} to show
that~(\ref{eq:est-proba}) holds true if $t-s>s^{(1)}_t$. So let us define
$$
s_t(b)=bt,\qquad\text{where}\qquad b\in(0,1/2).
$$
By Proposition~\ref{prop:a=t}, we immediately have $F(t,x_t(c),t-s_t(b))=o(1)$.
We observe that
$$
G(t-s_t(b),x_t(c))=W(t-s_t(b))\int_{s^{(1)}_t}^{t-s_t(b)}\frac{e^{-\alpha x}}{W_\alpha(x)}e^{(\lceil
  x_t(c)\rceil-1)\log(1-1/W_\alpha(x))}(\alpha dx+\delta_{t-s_t(b)}(dx)).
$$
Using the inequality $\log(1-x)\leq -x$, Lemma~\ref{lem:W_theta}~(iii) entails that the
contribution of the Dirac mass is
$$
O\left(\frac{1}{t}\,
e^{-x_t(c)/W_\alpha(t-s_t(b))}\right).
$$
Fix $\eta\in(0,1)$. Using the expression of $x_t(c)$ and the fact that
$1/W_\alpha(t)\geq(1-\eta)\psi'(\alpha)/(\alpha t)$ for $t$ large enough, this last
quantity is
$$
O\left(e^{-\frac{\alpha(1-\eta)}{4(1-b)}\,t}\,t^{-\left(1-\frac{1-\eta}{4(1-b)}\right)}\right)=o\left(e^{-\alpha
  s_t(b)}\right),
$$
where the last equality is valid if one chooses $\eta<(1-2b)^2$.

Hence
$$
G(t-s_t(b),x_t(c))=W(t-s_t(b))\int_{s^{(1)}_t}^{t-s_t(b)}\frac{e^{-\alpha x}}{W_\alpha(x)}e^{(\lceil
  x_t(c)\rceil-1)\log(1-1/W_\alpha(x))}\alpha dx+o\left(e^{-\alpha
  s_t(b)}\right).
$$
Now, the integral in the r.h.s.\ is exactly the same as in~(\ref{eq:pf-crit-lim-esp}), except
for the interval of integration. We actually proved in the proof of Proposition~\ref{prop:a=t}
that, since $(1-b)>1/2$, this integral is equivalent to
$$
\int_{\frac{1-\varepsilon}{2}t}^{\frac{1+\varepsilon}{2}t}\frac{e^{-\alpha x}}{W_\alpha(x)}e^{(\lceil
  x_t(c)\rceil-1)\log(1-1/W_\alpha(x))}\alpha dx,
$$
which is itself equivalent to
$$
\frac{1}{W(t)}\,\sqrt{\frac{2\pi}{\alpha}}\,e^{B-\frac{\psi'(\alpha)}{2}}\,e^{-\alpha c}.
$$
Therefore,
\begin{equation}
  \label{eq:pf-crit-freq-interm}
  G(t-s_t(b),x_t(c))\sim e^{-\alpha
    s_t(b)}\,\sqrt{\frac{2\pi}{\alpha}}\,e^{B-\frac{\psi'(\alpha)}{2}}\,e^{-\alpha c},
\end{equation}
and, recalling that~(\ref{eq:est-proba}) holds (with our current notation), the proof of
Lemma~\ref{lem:crit-freq} will be completed if we can prove that
$$
G(t-s_t(b),x_t(c))\sim F(t-s_t(b),x_t(c))
$$
as $t\rightarrow+\infty$. Again, this is implied by the estimate
\begin{equation}
  \label{eq:goal-goal-goal}
  \EE_t[M_{t-s_t(b)}(x_t(c),s^{(1)}_t)(M_{t-s_t(b)}(x_t(c),s^{(1)}_t)-1)]=o(e^{-\alpha s_t(b)}) ,
\end{equation}
which we now prove.

Applying Lemma~\ref{lem:summary-technique}~(\ref{eq:jolie-equation}) with $x=x_t(c)$,
$s_1=s^{(1)}_t$ and $s_2=t-s_t(b)$, and combining the result with~(\ref{eq:terme-qui-sauve})
and the fact that, for all $s\leq t$,
$$
e^{-\alpha s}+\int_0^s \alpha e^{-\alpha z}\frac{W_\alpha(t)-W_\alpha(s-z)}{W_\alpha(t)}dz\leq
e^{-\alpha s}+\int_0^s \alpha e^{-\alpha z}=1,
$$
we obtain
\begin{multline}
  \EE_t[M_{t-s_t(b)}(x_t(c),s^{(1)}_t)(M_{t-s_t(b)}(x_t(c),s^{(1)}_t)-1)] \\ \leq C\, \EE_t
  K_{t-s_t(b)}(x_t(c),s^{(1)}_t,t-s_t(b))\,\left(1+e^{-\alpha
      s^{(1)}_t}\,e^{\alpha(t-s_t(b))}\right)\times \\
  \left[(1+t)\left(1-\frac{1}{W_\alpha(t-s_t(b))}\right)^{\lceil
      x_t(c)\rceil-1}+e^{-\alpha s^{(1)}_t}\,e^{\alpha(t-s_t(b))}\,\EE_t
    K_{t-s_t(b)}\left(\frac{x_t(c)}{2},s^{(1)}_t,t-s_t(b)\right)\right].
  \label{eq:calcul-calcul-calcul}
\end{multline}
Fix again $\eta\in(0,1)$. Using the inequality $\log (1-x)\leq -x$ and
Lemma~\ref{lem:W_theta}~(iii), we have for $t$ large enough,
\begin{align*}
  \left(1-\frac{1}{W_\alpha(t-s_t(b))}\right)^{\lceil x_t(c)\rceil-1} &
  \leq C\exp\left(-\frac{\alpha^2}{4\psi'(\alpha)}\left(t^2-\frac{t\log
        t}{\alpha}+2ct\right)\frac{(1-\eta)\psi'(\alpha)}{\alpha(t-s_t(b))}\right) \\ &
  \leq C \exp\left(-\frac{\alpha(1-\eta)}{4(1-b)}\,t\right)\, t^{1/2},
\end{align*}
where we used the inequality $1/(1-b)<2$ to upper bound the exponent of $t$ in the last
inequality. Using Lemma~\ref{lem:summary-technique}~(\ref{eq:K-version-1}), this last
inequality yields
\begin{align*}
  \EE_t K_{t-s_t(b)}(x_t(c),s^{(1)}_t,t-s_t(b)) &
  \leq\frac{b}{\alpha}\int_{s^{(1)}_t}^{t-s_t(b)}\left(1-\frac{1}{W_\alpha(y)}\right)^{\lceil
    x_t(c)\rceil-1}(\alpha dy+\delta_0(dy)) \\
  & \leq\frac{b}{\alpha}\left(1-\frac{1}{W_\alpha(t-s_t(b))}\right)^{\lceil x_t(c)\rceil-1}
  \int_{s^{(1)}_t}^{t-s_t(b)}(\alpha dy+\delta_0(dy)) \\ & \leq
  C\exp\left(-\frac{\alpha(1-\eta)}{4(1-b)}\,t\right)\, t^{3/2}.
\end{align*}
Similarly,
\begin{align*}
  \EE_t K_{t-s_t(b)}(x_t(c)/2,s^{(1)}_t,t-s_t(b)) &
  \leq\frac{b}{\alpha}\left(1-\frac{1}{W_\alpha(t-s_t(b))}\right)^{\lceil x_t(c)/2\rceil-1}
  \int_{s^{(1)}_t}^{t-s_t(b)}(\alpha dy+\delta_0(dy)) \\
  & \leq C \exp\left(-\frac{\alpha(1-\eta)}{8(1-b)}\,t\right)\, t^{5/4}.
\end{align*}
Combining the previous inequalities with~(\ref{eq:calcul-calcul-calcul}), we finally obtain
\begin{multline*}
  \EE_t[M_{t-s_t(b)}(x_t(c),s^{(1)}_t)(M_{t-s_t(b)}(x_t(c),s^{(1)}_t)-1)] \\
  \leq C\, t^3\, e^{-\alpha s_t(b)} 
  \exp\left(-\alpha
    t\left(\frac{1-\varepsilon}{2}+\frac{1-\eta}{4(1-b)}-1\right)\right) \times \\
  \left[\exp\left(-\alpha t\,\frac{1-\eta}{4(1-b)}\right)+
    \exp\left(-\alpha
    t\left(b+\frac{1-\varepsilon}{2}+\frac{1-\eta}{8(1-b)}-1\right)\right)\right].
\end{multline*}
Remember now that $\varepsilon$ and $\eta$ are free parameters in $(0,1)$. We may assume that
they are linked to $b$ by the equation
$$
\frac{1-\eta}{4(1-b)}=\frac{1-\varepsilon}{2},\qquad\text{or}\qquad
b=\frac{1}{2}-\frac{\varepsilon-\eta}{2(1-\varepsilon)}
$$
which is always possible since $b<1/2$. This yields
\begin{multline*}
  \EE_t[M_{t-s_t(b)}(x_t(c),s^{(1)}_t)(M_{t-s_t(b)}(x_t(c),s^{(1)}_t)-1)]
  \\ \leq C\, t^3\, e^{-\alpha s_t(b)} \left[\exp\left(-\alpha
      t\,\frac{1-3\varepsilon}{2}\right)
    +\exp\left(-\alpha
      t\left(\frac{1}{4}-\frac{\varepsilon-\eta}{2(1-\varepsilon)}
        -\frac{7\varepsilon}{4}\right)\right)\right].
\end{multline*}
Taking $b$ close enough to $1/2$ allows to take both $\varepsilon$ and $\eta$ as close to 0 as
desired. Therefore,~(\ref{eq:goal-goal-goal}) is proved and the proof of
Lemma~\ref{lem:crit-freq} is completed.\hfill$\Box$
\bigskip

\subsection{Old haplotypes}
\label{sec:old-crit}

\begin{thm}
  \label{thm:cv-distr-old-a=t}
  Assume $\alpha=\theta$ and define for all $a\in\RR$
  $$
  x_t(a)=t-\frac{\log t}{\alpha}+a.
  $$
  For all $n\in\NN$, $k_1,\ldots,k_n\in\ZZ_+$ and $a_1>a_2>\ldots>a_n$, we have
  \begin{multline}
    \lim_{t\rightarrow+\infty}\,\PP_t\Big[O_t(x_t(a_1))=k_1,\ O_t(x_t(a_2))-O_t(x_t(a_1))=k_2,\ldots,\
    O_t(x_t(a_n))-O_t(x_t(a_{n-1}))=k_n\Big] \\
    =\alpha\,\binom{k_1+\ldots+k_n}{k_1,\ldots,k_n}\,
    \frac{e^{-\alpha k_1 a_1}(e^{-\alpha a_2}-e^{-\alpha a_1})^{k_2}\ldots(e^{-\alpha a_n}-e^{-\alpha a_{n-1}})^{k_n}}
    {(\alpha+e^{-\alpha a_n})^{k_1+\ldots+k_n+1}}. \label{eq:goal-old-cas-critique}
  \end{multline}
  In addition, the family of ${\cal M}(\RR)$-valued random variables $(Z_t,t\geq 0)$, defined for all $t\geq
  0$ by
  $$
  Z_t:=\sum_{k\geq 1}\delta_{A^{(k)}_t-t+\frac{\log t}{\alpha}},
  $$
  converges as $t\rightarrow+\infty$ in $\PP^\star$-distribution in ${\cal M}(\RR)$ equipped with the semi-vague topology to a mixed
  Poisson point measure on $\RR$ with intensity measure
  \begin{equation}
    \label{eq:intens-measure-P-age-cas-critique}
    {\cal E}\,e^{-\theta a}\,da,    
  \end{equation}
  where the mixture coefficient ${\cal E}$ has exponential distribution with parameter $1$.
\end{thm}

Again, this result follows from the next lemma exactly as Theorem~\ref{thm:cv-distr-old-a<t} followed from
Lemma~\ref{lem:age}.

\begin{lem}
  \label{lem:age-cas-critique}
  For all $a\in\RR$
  $$
  \lim_{t\rightarrow+\infty}\,\PP_t\left(A^{(1)}_t\leq t-\frac{\log
      t}{\alpha}+a\right)=\frac{1}{1+\frac{1}{\alpha}\,e^{-\alpha a}}.
  $$
\end{lem}

\paragraph{Proof.}

We define $F(t,x)$ and $G(t,x)$ exactly as in the proof of Lemma~\ref{lem:age} and we put
$$
s_t(b)=\frac{b}{\alpha}\log t\qquad\text{with}\qquad b\in(0,1).
$$
With this new notation,~(\ref{eq:chouette}) holds true and Proposition~\ref{prop:a=t-age} implies that
$G(t,t-s_t(b))=o(1)$. In addition, one checks exactly as in the proof of Proposition~\ref{prop:a=t-age} that
$$
G(t-s_t(b),x_t(a))\sim\frac{e^{-\alpha a}}{\alpha}\,e^{-\alpha s_t(b)}
$$
when $t\rightarrow+\infty$. The proof will then be completed if we can
prove~(\ref{eq:goal-calcul-age-loi}). We first observe that $W_\alpha'(x)=e^{-\alpha x}W'(x)$ is bounded
thanks to~(\ref{eq:link-W-excursion}). Therefore, by Lemma~\ref{lem:W_theta}~(iii),
\begin{multline*}
  e^{-\alpha x_t(a)}+\int_0^{x_t(a)}\alpha e^{-\alpha
    z}\frac{W_\alpha(t-s_t(b))-W_\alpha(x_t(a)-z)}{W_\alpha(t-s_t(b))}dz \\
  \begin{aligned}
    & \leq C\left(te^{-\alpha
        t}+\int_0^{x_t(a)}e^{-\alpha z}\frac{t-s_t(b)-x_t(a)+z}{t-s_t(b)}dz\right) \\ & \leq
    C\frac{\log t}{t}.
  \end{aligned}
\end{multline*}
Combining this inequality with~(\ref{eq:terme-qui-sauve}) and
Lemma~\ref{lem:summary-technique}~(\ref{eq:jolie-equation}) in which we take $x=1$,
$s_1=x_t(a)$ and $s_2=t-s_t(b)$ yields
\begin{multline}
  \EE_t\left[O_{t-s_t(b)}(x_t(a))\left(O_{t-s_t(b)}(x_t(a))-1\right)\right] \leq C
  \EE_t K_{t-s_t(b)}(1,x_t(a),t-s_t(b))\left(1+e^{\alpha(t-s_t(b)-x_t(a))}\right)\times \\
  \left[(1+\log t)\frac{\log t}{t}+e^{\alpha(t-s_t(b)-x_t(a))}\EE_t
    K_{t-s_t(b)}(1,x_t(a),t-s_t(b))\right]. \label{eq:calcul-age-loi-cas-critique}
\end{multline}
By Lemma~\ref{lem:summary-technique}~(\ref{eq:K}), we have
\begin{equation*}
  \EE_t K_{t-s_t(b)}(1,x_t(a),t-s(b))\leq\frac{b}{\alpha}\int_{x_t(a)}^{t-s_t(b)}(\alpha
  dy+\delta_{t-s_t(b)}(dy))\left(e^{-\alpha y}+\int_0^{y}\alpha e^{-\alpha z}
    \frac{W_\alpha(y)-W_\alpha(y-z)}{W_\alpha(y)}dz\right).
\end{equation*}
Using again the fact that $W_\alpha'(x)$ is bounded and that $1/W_\alpha(y)\leq C/y$ for all
$y$ large enough, we deduce that
\begin{align*}
  \EE_t K_{t-s_t(b)}(1,x_t(a),t-s_t(b)) & \leq
  \int_{x_t(a)}^{t-s_t(b)}(\alpha dy+\delta_{t-s_t(b)}(dy))\left(e^{-\alpha y}+\frac{C}{y}\int_0^{y}
    ze^{-\alpha z}dz\right) \\ & \leq
  C\int_{x_t(a)}^{t-s_t(b)}\frac{1}{y}(\alpha dy+\delta_{t-s_t(b)}(dy)) \\ & \leq C
  \frac{\log t}{t}.
\end{align*}
Therefore, it follows
from~(\ref{eq:calcul-age-loi-cas-critique}) that
$$
\EE_t\left[O_{t-s_t(b)}(x_t(a))\left(O_{t-s_t(b)}(x_t(a))-1\right)\right]\leq C\,e^{-\alpha
  s_t(b)}\,\log t\,\left(\frac{(\log t)^2}{t}+\frac{\log t}{t^b}\right)=o(e^{-\alpha s_t(b)}).
$$
This completes the proof of Lemma~\ref{lem:age-cas-critique}.\hfill$\Box$
\bigskip

\appendix

\section{Proof of Lemma~\ref{lem:summary-technique}}

Recall the notation introduced in Section~\ref{sec:prelim-lem}.

As seen in the proof of Theorem~\ref{thm:cv-loi-most-frequent}, this result (and actually all
the results of Sections~\ref{sec:subcrit} and~\ref{sec:crit}) are consequences of estimates of
the form
$$
\PP_t(M_t(x_t,s_t)\geq 1)\sim\EE_t M_t(x_t,s_t)
$$
as $t\rightarrow+\infty$, for convenient choices of $x_t$ and $s_t$. We chose to prove this
result using the inequality
$$
0\leq\EE_t M_t(x_t,s_t)-
\PP_t(M_t(x_t,s_t)\geq 1)\leq\EE_t[M_t(x_t,s_t)(M_t(x_t,s_t)-1)],
$$
i.e.\ proving that
$$
\EE_t[M_t(x_t,s_t)(M_t(x_t,s_t)-1)]=o\left(\EE_t M_t(x_t,s_t)\right).
$$
Such results are obtained using Lemma~\ref{lem:summary-technique}, which is an immediate
consequence of the following two lemmas.
\bigskip

We need to define the random variable $K'_t(x,s_1,s_2)$ by slightly modifying the definition
of $K_t(x,s_1,s_2)$: introducing an independent random variable $H'$ distributed as $H_i$
conditional on $\{H_i< t\}$, $K'_t(x,s_1,s_2)$ is the number of haplotypes carried by more
than $x$ individuals alive at time $t$, whose last mutation occurred on branch $0$, and is
older than $s_1$ and younger than $s_2\wedge H'$.

As a first step, we compute an upper bound of $\EE_t
[M_t(x,s_1,s_2)(M_t(x,s_1,s_2)-1)]$ expressed in terms of $K_t$ and $K'_t$.

\begin{lem}
  \label{lem:technique}
  For all $t>0$, $x\geq 1$, $0\leq s_1\leq s_2\leq+\infty$, we have
  \begin{align*}
    \EE_t [M_t(x,s_1,s_2)(M_t(x,s_1,s_2)-1)] & \le \EE_t[K_t(x,s_1,s_2)(K_t(x,s_1,s_2)-1)] \\ & +
    (\EE_t N_t)\EE_t[K'_t(x,s_1,s_2)(K'_t(x,s_1,s_2)-1)] \\ & +8 (\EE_t N_t)(\EE_t K_t(\lceil
    x/2\rceil,s_1,s_2))(\EE_t K'_t(x,s_1,s_2)) \\ & +8 (\EE_t N_t)^2(\EE_t K'_t(\lceil
    x/2\rceil,s_1,s_2))(\EE_t K'_t(x,s_1,s_2)).
  \end{align*}
\end{lem}

\paragraph{Proof.}
We let $M_i$ be the number of mutations on branch $i$ (this branch has length $H_i$),
considering only the mutations younger than $t$ when $i=0$. For all $j\le M_i$,
we define $\ell_{ij}$ the duration elapsed \emph{since} the $j$-th oldest mutation on branch $i$, with
$\ell_{i(M_i+1)}=0$, $\ell_{i0}=H_i$ and $\ell_{00}=t$. We also define $M'_0$ as the smallest $k\geq 1$ such
that $\ell_{0k}\leq H'$ (and $M'_0=0$ there is no such $k\geq 1$).

For $0\leq j\leq M_i$, denote by $R_t^{ij}$ the number of individuals alive at time $t$ descending clonally
from the time interval $I_{ij}:=(t-\ell_{ij}, t-\ell_{i(j+1)})$ on branch $i$. More specifically, for a
progenitor individual alive on the time interval $(a,b)$ and experiencing no mutation between times $a$ and
$b$, we refer to `clonal descendants from the time interval $(a,b)$' as those individuals alive at $t$
(including possibly the progenitor) descending from those daughters of the progenitor who were born during the
time interval $(a,b)$, and that still carry the same type the progenitor carried at time $a$. Using the
notation $A_{ij}=A_{ij}(t,x,s_1,s_2):=\{R_t^{ij}\geq x,\,\ell_{ij}\in[s_1,s_2)\}$, we have
\begin{gather*}
  M_t(x,s_1,s_2)=\sum_{0\leq j\leq M_0}\mathbbm{1}_{A_{0j}}+\sum_{1\le i< N_t}\sum_{1\le j\le
    M_i}\mathbbm{1}_{A_{ij}}, \\
  K_t(x,s_1,s_2)=\sum_{0\leq j\leq M_0}\mathbbm{1}_{A_{0j}}
\end{gather*}
and
$$
K'_t(x,s_1,s_2)=\sum_{0\leq j\leq M'_0}\mathbbm{1}_{A_{0j}}.
$$

Therefore, by construction of the coalescent point process,
\begin{align*}
  \EE_t M_t(x,s_1,s_2) & =\EE_t K_t(x,s_1,s_2)+\sum_{i\geq 1}\sum_{j\geq 1}\PP_t(A_{ij},i<N_t,j\leq
  M_i) \\ & =\EE_t K_t(x,s_1,s_2)+\sum_{i\geq 1}\sum_{j\geq 1}\PP_t(A_{ij},j\leq
  M_i\mid i<N_t)\PP_t(i<N_t) \\ & =\EE_t K_t(x,s_1,s_2)+\sum_{i\geq 1}\sum_{j\geq 1}\PP_t(A_{0j},j\leq
  M'_0)\PP_t(i<N_t) \\ & =\EE_t K_t(x,s_1,s_2)+(\EE_t N_t-1)\EE_t K'_t(x,s_1,s_2).
\end{align*}
Now,
\begin{align*}
  M_{t}(x,s_1,s_2)(M_{t}(x,s_1,s_2)-1) & =2\sum_{0\leq k< j\leq M_0}\mathbbm{1}_{A_{0k}}\mathbbm{1}_{A_{0j}}
  +2\sum_{i=1}^{N_t-1}\sum_{1\leq k<j\leq M_i}\mathbbm{1}_{A_{ik}}\mathbbm{1}_{A_{ij}} \\ &
  +2\sum_{k=0}^{M_0}\sum_{i=1}^{N_t-1}\sum_{j=1}^{M_i}\mathbbm{1}_{A_{0k}}\mathbbm{1}_{A_{ij}}
  +2\sum_{1\leq l<i< N_t}\sum_{k=1}^{M_l}\sum_{j=1}^{M_i}\mathbbm{1}_{A_{lk}}\mathbbm{1}_{A_{ij}}.
\end{align*}
Hence, using a similar computation as above,
\begin{multline*}
  \EE_t M_{t}(x,s_1,s_2)(M_{t}(x,s_1,s_2)-1) \\
  \begin{aligned}
    & = \EE_t K_t(x,s_1,s_2)(K_t(x,s_1,s_2)-1)+
    2\sum_{i\geq 1}\sum_{k\geq 1}\sum_{j>k}\PP_t(A_{0k}\cap A_{0j},j\leq M'_0)\PP_t(i<N_t) \\
    & +2\sum_{k\geq 0}\sum_{i\geq 1}\sum_{j\geq 1}\PP_t(A_{0k}\cap A_{ij}, k\leq M_0,i< N_t,j\leq M_i) \\
    & +2\sum_{l\geq 1}\sum_{k\geq 1}\sum_{i>l}\sum_{j\geq 1}\PP_t(A_{0k}\cap A_{(i-l)j},k\leq
    M'_0,i-l<N_t,j\leq M_{i-l})\PP_t(l<N_t) \\
    & =\EE_t K_t(x,s_1,s_2)(K_t(x,s_1,s_2)-1)+(\EE_t N_t-1)\EE_t
    K'_t(x,s_1,s_2)(K'_t(x,s_1,s_2)-1) \\
    & +2\sum_{k\geq 0}\sum_{i\geq 1}\sum_{j\geq 1}\PP_t(A_{0k}\cap A_{ij}, k\leq M_0,i< N_t,j\leq M_i) \\
    & +2(\EE_t N_t-1)\sum_{k\geq 1}\sum_{i\geq 1}\sum_{j\geq 1}\PP_t(A_{0k}\cap A_{ij}, k\leq M'_0,i<
    N_t,j\leq M_i).
  \end{aligned}
\end{multline*}
For short, we write 
\begin{multline}
  \EE_t M_{t}(x,s_1,s_2)(M_{t}(x,s_1,s_2)-1) \\
  \leq\EE_t K_t(x,s_1,s_2)(K_t(x,s_1,s_2)-1)+(\EE_t N_t)\EE_t
  K'_t(x,s_1,s_2)(K'_t(x,s_1,s_2)-1) \\
  +2\sum_{k\geq 0}\sum_{i\geq 1}\sum_{j\geq 1}\PP_t(A_{0k}\cap A_{ij}\cap B_{ijk})
  +2(\EE_t N_t)\sum_{k\geq 1}\sum_{i\geq 1}\sum_{j\geq 1}\PP_t(A_{0k}\cap A_{ij}\cap B'_{ijk})
  \label{eqn : esperance L(L-1)}
\end{multline}
where $B_{ijk}:=\{ k\le M_0, i<N_{t}, j \le M_{i}\}$ and $B'_{ijk}:=\{ k\le M'_0, i<N_{t}, j \le M_{i}\}$.

Now for any positive integers $i,k$, define the three following events
$$
\alpha_{ik}:=\{\max_{1\le j\le i} H_j >\ell_{0k}\},\qquad \beta_{ik}:=\{\ell_{0(k+1)}<\max_{1\le j\le i} H_j  \le \ell_{0k}\},\qquad \gamma_{ik}:=\{\max_{1\le j\le i} H_j  \le \ell_{0(k+1)}\}.
$$ 
We are going to state and prove six inequalities, where the left-hand side is obtained by
intersecting each event $A_{0k}\cap A_{ij}\cap B_{ijk}$ or $A_{0k}\cap A_{ij}\cap B'_{ijk}$
with each of the preceding ones $\alpha_{ik}$, $\beta_{ik}$, $\gamma_{ik}$, and summing over
$i,j\geq 1$ and $k\geq 0$ for the events involving $B_{ijk}$, and $k\geq 1$ for the events
involving $B'_{ijk}$.
\begin{align}
  \sum_{i,j,k}\PP_t (\alpha_{ik}\cap A_{0k}\cap A_{ij}\cap B_{ijk}) & \le (\EE_t N_t)(\EE_t
  K_{t}(x,s_1,s_2))(\EE_t K'_t(x,s_1,s_2)), \label{eqn : alpha} \\
  \sum_{i,j,k}\PP_t (\alpha_{ik}\cap A_{0k}\cap A_{ij}\cap B'_{ijk}) & \le (\EE_t N_t)(\EE_t
  K'_t(x,s_1,s_2))^2, \label{eqn : alpha'} \\
  \sum_{i,j,k}\PP_t (\gamma_{ik}\cap A_{0k}\cap A_{ij}\cap B_{ijk}) & \le (\EE_t N_t)(\EE_t
  K_{t}(x,s_1,s_2))(\EE_t K'_t(x,s_1,s_2)), \label{eqn : gamma} \\
  \sum_{i,j,k}\PP_t (\gamma_{ik}\cap A_{0k}\cap A_{ij}\cap B'_{ijk}) & \le (\EE_t N_t)(\EE_t
  K'_t(x,s_1,s_2))^2, \label{eqn : gamma'} \\
  \sum_{i,j,k}\PP_t (\beta_{ik}\cap A_{0k}\cap A_{ij}\cap B_{ijk}) & \le 2(\EE_t N_t)(\EE_t
  K_t(\lceil x/2\rceil,s_1,s_2))(\EE_t K'_t(x,s_1,s_2)), \label{eqn : beta} \\
  \sum_{i,j,k}\PP_t (\beta_{ik}\cap A_{0k}\cap A_{ij}\cap B'_{ijk}) & \le 2(\EE_t N_t)(\EE_t
  K'_t(\lceil x/2\rceil,s_1,s_2))(\EE_t K'_t(x,s_1,s_2)). \label{eqn : beta'}
\end{align}
Combining these six equations with \eqref{eqn : esperance L(L-1)} and with the inequalities
$K_t(x,s_1,s_2)\leq K_t(\lceil x/2\rceil,s_1,s_2)$ and $K'_t(x,s_1,s_2)\leq K'_t(\lceil x/2\rceil,s_1,s_2)$
yields the inequality given in the lemma.

We are going to detail the proof of the inequalities \eqref{eqn : alpha}, \eqref{eqn : gamma} and
\eqref{eqn : beta} (in this order). The other inequalities can be proved using the same computations.
Let us start with \eqref{eqn : alpha}.  Hereafter we denote by $A_{0k}^{(i)}$ the event $\{I_{0k} \mbox{ has
  more than }x \mbox{ clonal descendants within }\{0, \ldots, i-1\}\text{\ and\ }\ell_{0k}\in[s_1,s_2)\}$.
\begin{align}
  \PP_t (\alpha_{ik}\cap A_{0k}\cap A_{ij}\cap B_{ijk}) & = \PP_t (\alpha_{ik}\cap A_{0k}^{(i)}\cap
  A_{ij}\cap B_{ijk}) \notag \\
  & \le \PP_t (k\le M_0,A_{0k}^{(i)}, i<N_{t}, j \le M_{i}, A_{ij} ) \notag \\
  & = \PP_t (k\le M_0,A_{0k}^{(i)}, i<N_{t}, j \le M_{i}, A_{ij} \mid i< N_t)\PP_t(i< N_t) \notag \\
  & =\PP_t (k\le M_0, A_{0k}^{(i)}\mid i<N_{t}) \PP_t( A_{ij}, j\le M_i \mid i< N_{t})\PP_t (i<
  N_t) \notag \\
  & =\PP_t (k\le M_0, A_{0k}^{(i)}, i<N_{t}) \PP_t( A_{0j}, j\le M'_0). \label{eq:tech-1}
\end{align}
Then, denoting $\rho_k$ the index of the $\lceil x \rceil$-th individual carrying the type of
interval $I_{0k}$ ($:=+\infty$ if there is no such individual),
\begin{align*}
  \sum_{i,j,k}\PP_t (\alpha_{ik}\cap A_{0k}\cap A_{ij}\cap B_{ijk}) & \le \EE_t \left(\sum_{0\le k\le M_0} \sum_{i\ge 1}\mathbbm{1}_{\ell_{0k}\in[s_1,s_2)}\mathbbm{1}_{\rho_k\le i<N_{t}}\right) \EE_t \sum_{j\le M'_0}\mathbbm{1}_{A_{0j}}\\
  & =\EE_t \left(\sum_{0\le k\le M_0} \mathbbm{1}_{\ell_{0k}\in[s_1,s_2)}(N_{t}-\rho_k)^+\right) \EE_t K'_{t}(x,s_1,s_2). 
\end{align*}
Now, the lack of memory property of geometric distributions yields
\begin{multline*}
  \EE_t \left(\sum_{0\le k\le M_0}\mathbbm{1}_{\ell_{0k}\in[s_1,s_2)} (N_{t}-\rho_k)^+\right) \\
  \begin{aligned}
    & =\sum_{k\ge 0} \EE_t\left( N_{t}-\rho_k \mid k\le M_0, \ell_{0k}\in[s_1,s_2),\rho_k <\infty \right) \PP_t\left( k\le M_0,\ell_{0k}\in[s_1,s_2), \rho_k <\infty \right)\\
    & =(\EE_t N_{t}) \sum_{k\ge 0}\PP_t\left( k\le M_0, A_{0k} \right)\\
    & =(\EE_t N_{t}) (\EE_t K_{t}(x,s_1,s_2)),
  \end{aligned}
\end{multline*}
which entails~\eqref{eqn : alpha}.

Next, let us proceed with \eqref{eqn : gamma} and let $\sigma_k$ denote the label of the first branch with length greater than $\ell_{0(k+1)}$. Observe that conditionally on $\ell_{0(k+1)}$, $A_{0k}$ is independent of the branch lengths occurring before $\sigma_k$, and further, the events $\{i\le\sigma_k\}=\{\max_{j\le i-1} H_j \le\ell_{0(k+1)}\}$, $\{H_i\le\ell_{0(k+1)}, j\le M_i , A_{ij} \}$ and  $\{k\le M_0 ,A_{0k} \}$ are independent. As a consequence,
\begin{multline*}
  \PP_t (\gamma_{ik}\cap A_{0k}\cap A_{ij}\cap B_{ijk}) = \PP_t (i<\sigma_k, k\le M_0, A_{0k}, j\le
  M_i, A_{ij})\\
  \begin{aligned}
    & = \PP_t (i\le\sigma_k, k\le M_0, A_{0k}, H_i\le\ell_{0(k+1)} , j\le M_i, A_{ij})\\
    & =	\EE_t (\PP_t(i\le\sigma_k\mid\ell_{0(k+1)}) \PP_t(A_{0k},k\le M_0\mid \ell_{0(k+1)}) \PP_t(H_i\le\ell_{0(k+1)} , j\le M_i, A_{ij} \mid \ell_{0(k+1)}))\\
    & \leq\EE_t (\PP_t(i\le\sigma_k, A_{0k},k\le M_0\mid \ell_{0(k+1)}) \PP_t(j\le M_i, A_{ij} \mid \ell_{0(k+1)}))\\
    & =\EE_t (\PP_t(i\le\sigma_k, A_{0k},k\le M_0\mid \ell_{0(k+1)}) \PP_t( j\le M_i, A_{ij}))\\
    & =	\PP_t(j\leq M'_0,A_{0j})\PP_t (i\le\sigma_k, A_{0k},k\le M_0).
  \end{aligned}
\end{multline*}
As a consequence, since $\sigma_k$ and $\{k\leq M_0,\,A_{0k}\}$ are independent conditionally on $\ell_{0(k+1)}$,
\begin{align*}
  \sum_{i,j,k}\PP_t (\gamma_{ik}\cap A_{0k}\cap A_{ij}\cap B_{ijk}) & \le (\EE_t K'_{t}(x,s_1,s_2))\sum_{ k\ge 0} \EE_t \left(\sum_{i\le \sigma_k } \mathbbm{1}_{k\le M_0,  A_{0k}}\right) \\
  & = (\EE_t K'_{t}(x,s_1,s_2))\sum_{ k\ge 0} \EE_t \left(\PP_t (k\le M_0, A_{0k}\mid\ell_{0(k+1)}) \EE_t(\sigma_k\mid \ell_{0(k+1)})\right) \\
  & \le (\EE_t K'_{t}(x,s_1,s_2))\sum_{ k\ge 0} \EE_t \left(\PP_t (k\le M_0, A_{0k}\mid\ell_{0(k+1)}) \EE_t(N_{t})\right)\\
  & =(\EE_t K'_{t}(x,s_1,s_2))(\EE_t K_t(x,s_1,s_2))(\EE_t N_{t}), 
\end{align*}
which is \eqref{eqn : gamma}. 

Finally, let us turn to \eqref{eqn : beta}. Denote by $A^{\prime(i)}_{0k}$ (resp.\ $A^{\prime\prime(i)}_{0k}$)
the event that there exists at least $\lceil x/2\rceil$ individual with label \emph{smaller} (resp.\
\emph{larger}) than $i$ descending clonally from the time interval $I_{0k}$ and that $\ell_{0k}\in[s_1,s_2)$.
Then
\begin{equation*}
  \PP_t (\beta_{ik}\cap A_{0k}\cap A_{ij}\cap B_{ijk})
  \le\PP_t (\beta_{ik}\cap A^{\prime(i)}_{0k}\cap A_{ij}\cap B_{ijk})
  +\PP_t (\beta_{ik}\cap A^{\prime\prime(i)}_{0k}\cap A_{ij}\cap B_{ijk}).
\end{equation*}
Let us deal with the first term of the right-hand side of this last inequality. Exactly as in the proof of \eqref{eqn : alpha}, 
\begin{align*}
  \PP_t (\beta_{ik},k\le M_0, A^{\prime(i)}_{0k},  j\le M_i, A_{ij}, i< N_{t}) & \le\PP_t (k\le M_0,
  A^{\prime(i)}_{0k}, j\le M_i, A_{ij}, i< N_{t}) \\
  & =\PP_t (k\le M_0, A^{\prime(i)}_{0k}, i<N_{t}) \PP_t( A_{ij}, j\le M_i \mid i< N_{t})\\
  & =\PP_t (k\le M_0, A^{\prime(i)}_{0k}, i<N_{t}) \PP_t( A_{0j}, j\le M'_0),
\end{align*}
and we finally get
\begin{equation}
  \label{eqn : beta 1}
  \sum_{i,j,k}\PP_t (\beta_{ik}\cap A^{\prime(i)}_{0k}\cap A_{ij}\cap B_{ijk})\le (\EE_t
  N_{t})(\EE_t K_{t}(\lceil x/2\rceil,s_1,s_2))(\EE_t K'_{t}(x,s_1,s_2)).
\end{equation}
As for the second term, we need to define $J_i$ the unique integer satisfying $\ell_{0(J_i+1)}<\max_{1\le j\le
  i} H_j  \le \ell_{0J_i}$ ($J_i:=+\infty$ on $\{i\ge N_{t}\}$ and $J_i=k$ on $\beta_{ik}$). Then
\begin{equation}
  \label{eq:calc-lem-prelim}
  \sum_{i,j,k}\PP_t (\beta_{ik}\cap A^{\prime\prime(i)}_{0k}\cap A_{ij}\cap B_{ijk})
  =\sum_{i,j}\PP_t ( A^{\prime\prime(i)}_{0J_i}, A_{ij},j\le M_i, i<N_{t}).
\end{equation}
Set also $\ell^*_i:=\ell_{iM_i}$ the age of the oldest mutation on branch $H_i$ ($\ell^*_i=0$ if
$M_i=0$). Then conditional on $\{i< N_{t}\}$ and on the value of $\ell^*_i$, the numbers of clonal descendants
$R_t^{ij}$ of the interval $I_{ij}$ and the number, say $K^{(i)}$, of haplotypes whose last mutation is older
than $\ell^*_i$ and $s_1$, younger than $s_2$, and occurred on lineage $0$, and with more than $\lceil
x/2\rceil$ alive clonal descendants with labels larger than $i$, are independent, so that
\begin{align}
  \PP_t ( A^{\prime\prime(i)}_{0J_i}, A_{ij},j\le M_i, i<N_{t}) & \le\PP_t ( K^{(i)} \ge 1, A_{ij},j\le M_i,
  i<N_{t}) \notag \\
  & =\EE_t(\mathbbm{1}_{i<N_{t}}\PP_t(K^{(i)} \ge 1\mid i <N_{t},\ell^*_i)\PP_t( A_{ij},j\le M_i
  \mid i<N_{t}, \ell^*_i)) \notag \\
  & \le\EE_t(\mathbbm{1}_{i<N_{t}}\PP_t(K_{t}(\lceil x/2\rceil,s_1,s_2)\ge 1)\PP_t( A_{ij},j\le M_i \mid i<N_{t},
  \ell^*_i)) \notag \\
  & =\PP_t (K_{t}(\lceil x/2\rceil,s_1,s_2)\geq 1) \PP_t(A_{ij},j\le M_i \mid
  i<N_{t})\PP_t(i<N_t) \\ &
  \leq (\EE_t K_t(\lceil x/2\rceil,s_1,s_2))\PP_t (A_{0j},j\leq M'_0)\PP_t(i<N_t). \label{eq:tech-2}
\end{align}
We finally obtain
\begin{equation*}
  \sum_{i,j,k}\PP_t (\beta_{ik}\cap A^{\prime\prime(i)}_{0k}\cap A_{ij}\cap B_{ijk})\le
  (\EE_t N_t)(\EE_t K_{t}(\lceil x/2\rceil,s_1,s_2))(\EE_t K'_t(x,s_1,s_2)),
\end{equation*}
which completes the proof of \eqref{eqn : beta} by summing the last inequality and
inequality~\eqref{eqn : beta 1}. 

The proof of~(\ref{eqn : beta'}) is very similar, but needs further explanation. Let us define
the events $A^{\prime(i)}_{0k}$ and $A^{\prime\prime(i)}_{0k}$ similarly as above, with the
additional condition that $\ell_{0k}\leq H'$. Then, we first prove that
\begin{equation*}
  \sum_{i,j,k}\PP_t (\beta_{ik}\cap A^{\prime(i)}_{0k}\cap A_{ij}\cap B'_{ijk})\le (\EE_t
  N_{t})(\EE_t K'_{t}(\lceil x/2\rceil,s_1,s_2))(\EE_t K'_{t}(x,s_1,s_2))
\end{equation*}
following the very same computation as for~(\ref{eqn : beta 1}). Next, we observe that
$\PP_t(A^{\prime\prime(i)}_{00})=0$ since $H'<t$ a.s.\ and $\ell_{00}=t$.
Therefore,~(\ref{eq:calc-lem-prelim}) also holds true with our new definition of
$A^{\prime\prime(i)}_{0k}$. Thus, defining $K^{(i)}$ as the number of haplotypes whose last
mutation is older than $\ell^*_i$ and $s_1$, younger than $s_2$ \textbf{and} $H'$, and
occurred on lineage $0$, and with more than $\lceil x/2\rceil$ alive clonal descendants with
labels larger than $i$, the computation of~(\ref{eq:tech-2}) is true, provided that
$K_t(\lceil x/2\rceil,s_1,s_2)$ is replaced by $K'_t(\lceil x/2\rceil,s_1,s_2)$. We then
obtain
\begin{equation*}
  \sum_{i,j,k}\PP_t (\beta_{ik}\cap A^{\prime\prime(i)}_{0k}\cap A_{ij}\cap B'_{ijk})\le (\EE_t
  N_{t})(\EE_t K'_{t}(\lceil x/2\rceil,s_1,s_2))(\EE_t K'_{t}(x,s_1,s_2)),
\end{equation*}
and the proof of~(\ref{eqn : beta'}) is completed.\hfill$\Box$
\bigskip

Lemma~\ref{lem:summary-technique} follows from the combination of the previous lemma with the
following estimates on $K'_t(x,s_1,s_2)$ and $K_t(x,s_1,s_2)$.

\begin{lem}
  \label{lem:K-K'}
  For all $t>0$, $x\geq 1$, $0\leq s_1<s_2\leq t$, we have
  \begin{equation}
    \label{eq:K'}
    \EE_t K'_t(x,s_1,s_2)\leq\frac{\frac{1}{W(s_1)}-\frac{1}{W(t)}}{1-\frac{1}{W(t)}}\
    \EE_t K_t(x,s_1,s_2),
  \end{equation}
  \begin{equation}
    \label{eq:K-lem}
    \EE_t K_t(x,s_1,s_2)
    \leq\frac{b}{\alpha}\,\int_{s_1}^{s_2}\left(1-\frac{1}{W_\theta(y)}\right)^{\lceil
      x\rceil-1}\left(e^{-\theta y}+\int_0^{y}\theta
      e^{-\theta z}\frac{W_\theta(y)-W_\theta(y-z)}{W_\theta(y)}dz\right)(\theta
    dy+\delta_t(dy)),
  \end{equation}
  \begin{equation}
    \label{eq:K'(K'-1)}
    \EE_t K'_t(x,s_1,s_2)(K'_t(x,s_1,s_2)-1)\leq
    \frac{\frac{1}{W(s_1)}-\frac{1}{W(t)}}{1-\frac{1}{W(t)}}\
    \EE_t K_t(x,s_1,s_2)(K_t(x,s_1,s_2)-1),
  \end{equation}
  \begin{multline}
    \EE_t K_t(x,s_1,s_2)(K_t(x,s_1,s_2)-1)\leq \frac{2b}{\alpha}\,(\EE_t
    K_t(x,s_1,s_2))\,(1+\theta(s_2-s_1))\times \\ \times\left(1-\frac{1}{W_\theta(s_2)}\right)^{\lceil x\rceil -1}\left(e^{-\theta s_1}+\int_0^{s_1}\theta e^{-\theta
        z}\frac{W_\theta(s_2)-W_\theta(s_1-z)}{W_\theta(s_2)}dz\right). \label{eq:K(K-1)}
  \end{multline}
\end{lem}

\paragraph{Proof.}

With the notation of the proof of Lemma~\ref{lem:technique}, we have
\begin{align*}
  \EE_t K'_t(x,s_1,s_2) & =\sum_{k\geq 1}\PP_t(A_{0k},k\leq M'_0) \\
  & =\sum_{k\geq 1}\PP_t(A_{0k},k\leq M_0,H'\geq\ell_{0k}) \\ & \leq\sum_{k\geq 1}\PP_t(A_{0k},k\leq
  M_0)\PP(H'\geq s_1) \\ & \leq \EE_t K_t(x,s_1,s_2) \PP(H\geq s_1\mid H<t),
\end{align*}
which is inequality~(\ref{eq:K'}). Similarly,
\begin{align*}
  \EE_t K'_t(x,s_1,s_2)(K'_t(x,s_1,s_2)-1) & =2\sum_{k\geq 1}\sum_{j>k}\PP_t(A_{0k},A_{0j},j\leq M'_0) \\
  & \leq\sum_{k\geq 1}\sum_{j>k}\PP_t(A_{0k},A_{0j},j\leq M_0)\PP(H'\geq s_1) \\ & \leq \EE_t
  [K_t(x,s_1,s_2)(K_t(x,s_1,s_2)-1] \PP(H\geq s_1\mid H<t),
\end{align*}
which is inequality~(\ref{eq:K'(K'-1)}).

For the two other inequalities, let us define $R_t^{(a,b)}$ the number of individuals alive at
time $t$ descending clonally from the time interval $(a,b)$. More specifically, given a
progenitor individual alive on the time interval $(a,b)$ and experiencing no mutation between
times $a$ and $b$, $R_t^{(a,b)}$ is the number of individuals alive at time $t$ (including
this progenitor if $b\geq t$) descending from those daughters of the progenitor who
were born during the time interval $(a,b)$, and that still carry the same type that the
progenitor carried at time $a$. Since $W_\theta$ is the scale function associated with the
clonal reproduction process, for all $k\geq 0$,
\begin{align}
  \PP(R_t^{(a,b)}=k) & =\PP(N^\theta_{t-a}=k\mid\zeta=b-a) \notag \\
  & =\PP(N^\theta_{t-a}\not=0\mid\zeta=b-a)\PP(N^\theta_{t-a}=k\mid N^\theta_{t-a}\not=0)
  \notag \\
  & =\left(1-\mathbbm{1}_{t>b}\frac{W_\theta(t-b)}{W_\theta(t-a)}\right)
  \left(1-\frac{1}{W_\theta(t-a)}\right)^{k-1}\frac{1}{W_\theta(t-a)},
  \label{eq:loi-de-R}
\end{align}
where $N^\theta$ is the population size process of a clonal splitting tree and $\zeta$ is the
lifetime of the progenitor. (This result is actually Eq.~(5.3) of~\cite{CL1}.)

Note that, by construction of the splitting tree, replacing in the definition of $R_t^{(a,b)}$
the progenitor individual alive on the time interval $(a,b)$ by a \emph{clonal lineage} alive
on the time interval $(a,b)$, does not change anything to the distribution of $R_t^{(a,b)}$.
By \emph{lineage} alive on the time interval $(a,b)$, we mean here a finite sequence of
individuals $(i_k)_{1\leq k\leq K}$ such that individual $i_1$ was alive at time $a$, individual
$i_K$ was alive at time $b$, and for all $1\leq k\leq K-1$, individual $i_{k+1}$ was born from
individual $i_k$ at some time $a_k$ such that $a_1>a$ and $a_{K-1}<b$. By \emph{clonal}
lineage alive on the time interval $(a,b)$ we mean in addition that for all $1\leq k\leq K$,
individual $i_k$ experienced no mutation during the time interval $(a_{k-1},a_k)$, where
$a_0=a$ and $a_K=b$. 

Now, by definition of $K_t(x,s_1,s_2)$, we have
$$
\EE_t K_t(x,s_1,s_2)=\sum_{k\geq
  0}\EE_t\left[\mathbbm{1}_{\ell_{0k}\in[s_1,s_2]}\,
  \PP_t\left(R_t^{0k}\geq x\mid \ell_{0j}, j\geq
    0\right)\right],
$$
where
\begin{align*}
  \PP_t\left(R_t^{0k}\geq x\mid \ell_{0j},\, j\geq
    0\right) & =\frac{\PP\left(R^{0k}_t\geq x,\ N_t\geq 1\mid \ell_{0k},\,
      \ell_{0(k+1)}\right)}{\PP(N_t\geq 1)} \\ & \leq \frac{b}{\alpha}\,\PP\left(R^{0k}_t\geq
    x,\ N_t\geq 1\mid \ell_{0k},\,
    \ell_{0(k+1)}\right),
\end{align*}
since $\PP(N_t\geq 1)\geq\PP(N_s\geq 1,\ \forall s\geq 0)$ and the survival probability of the
splitting tree is $\alpha/b$. Now, the event $\{N_t\geq 1,\ R^{0k}_t\geq x\}$ is the event
where $N_t\geq 1$ and the clonal lineage on branch $0$ on the time interval
$(t-\ell_{0k},t-\ell_{0(k+1)})$ has more than $x$ clonal descendants alive at time $t$.
Therefore,
\begin{equation}
  \label{eq:super-inegalite}
  \PP_t\left(R_t^{0k}\geq x\mid \ell_{0j},\, j\geq
    0\right)\leq \frac{b}{\alpha}\,\PP\left(R_t^{(t-\ell_{0k},t-\ell_{0(k+1)})}\geq x\mid
    \ell_{0k},\, \ell_{0(k+1)}\right).
\end{equation}
Now, for all $k\geq 0$, $t-\ell_{0k}$ is distributed as the minimum of $t$ and a sum of $k$
i.i.d.\ exponential random variables of parameter $\theta$, and $t-\ell_{0(k+1)}$ as the
minimum of $t$ and the sum of $t-\ell_{0k}$ and an exponential random variable of parameter
$\theta$, independent of $t-\ell_{0k}$. Therefore, it follows from~(\ref{eq:loi-de-R}) that
\begin{align*}
  \frac{\alpha}{b}\,\EE_t K_t(x,s_1,s_2)
  & \leq\int_0^\infty dz\,\theta e^{-\theta z}\,\mathbbm{1}_{s_2=t}\,\PP(R_t^{(0,z)}\geq x)
  \\ & \qquad\qquad
  +\sum_{k\geq 1}\int_0^\infty dz\,\theta e^{-\theta z}\int_0^z dy\,\frac{\theta^k
    y^{k-1}}{(k-1)!}\,\mathbbm{1}_{y\in[t-s_2,t-s_1]}\,\PP(R_t^{(y,z)}\geq x) \\ &
  =\mathbbm{1}_{s_2=t}\int_0^\infty dz\,\theta e^{-\theta
    z}\left(1-\mathbbm{1}_{z<t}\frac{W_\theta(t-z)}{W_\theta(t)}\right)
  \left(1-\frac{1}{W_\theta(t)}\right)^{\lceil x\rceil-1} \\
  & \qquad\qquad +\int_{t-s_2}^{t-s_1} dy\,\theta e^{\theta y}\int_y^\infty dz\,\theta e^{-\theta z}
  \left(1-\mathbbm{1}_{z<t}\frac{W_\theta(t-z)}{W_\theta(t-y)}\right)\left(1-\frac{1}{W_\theta(t-y)}\right)^{\lceil
    x\rceil-1} \\ & =\int_{t-s_2}^{t-s_1}(\theta
  dy+\delta_0(dy))\,\left(1-\frac{1}{W_\theta(t-y)}\right)^{\lceil x\rceil-1}\left(1-\int_y^t dz\,\theta
    e^{-\theta (z-y)}\frac{W_\theta(t-z)}{W_\theta(t-y)}\right).
\end{align*}
Equation~(\ref{eq:K-lem}) then follows from the changes of variables $z'=z-y$ and $y'=t-y$,
and the identity $1=e^{-\theta y}+\int_0^y \theta e^{-\theta z}dz$.

Finally, let us turn to~(\ref{eq:K(K-1)}): first,
\begin{equation}
  \label{eq:youpla-1}
  \EE_tK_t(x,s_1,s_2)(K_t(x,s_1,s_2)-1)=2\sum_{0\leq j<k}\PP_t (A_{0j}, A_{0k}, k\leq M_0).
\end{equation}
Now, fix $k>l\geq 0$. Since $\ell_{0j}>\ell_{0(j+1)}\geq \ell_{0k}>\ell_{0(k+1)}$, on the
event $\{A_{0j},A_{0k},k\leq M_0\}$, we have $\ell_{0(j+1)}-\ell_{0k}\leq s_2-s_1$,
$\ell_{0k}\leq s_2$ and $\ell_{0(k+1)}\geq s_1-(\ell_{0k}-\ell_{0(k+1)})$. Therefore,
using~(\ref{eq:super-inegalite}) as before,
\begin{align*}
  & \PP_t (A_{0j}, A_{0k}, k\leq M_0)
  = \EE_t\left[\mathbbm{1}_{\ell_{0j}\in[s_1,s_2]}\PP_t\left(R_t^{0j}\geq x\mid \ell_{0j},\,\ell_{0(j+1)}\right)
    \mathbbm{1}_{\ell_{0k}\in[s_1,s_2]}\PP_t\left(R_t^{0k}\geq x\mid \ell_{0k},\,\ell_{0(k+1)}\right)
  \right] \\ &
  \leq\EE_t\Bigg[\mathbbm{1}_{\ell_{0j}\in[s_1,s_2]}\,\PP_t\left(R_t^{0j}\geq x\mid\ell_{0j},\,\ell_{0(j+1)}\right)
    \mathbbm{1}_{\ell_{0(j+1)}-\ell_{0k}\leq s_2-s_1}\times \\ &
  \qquad\qquad\qquad \times \left.\mathbbm{1}_{\ell_{0k}\geq
      s_1}\,\frac{b}{\alpha}\,\left(1-\mathbbm{1}_{\ell_{0k}-\ell_{0(k+1)}<s_1}
      \frac{W_\theta(s_1-(\ell_{0k}-\ell_{0(k+1)}))}{W_\theta(s_2)}\right)
    \left(1-\frac{1}{W_\theta(s_2)}\right)^{\lceil x\rceil-1}\mathbbm{1}_{\ell_{0k}>0}\right],
\end{align*}
where the last indicator comes from the fact that $1-1/W_\theta(\ell_{0k})=0$ when
$\ell_{0k}=0$. Now, on the event $\{\ell_{0k}>0\}$, one has
$$
\ell_{0n}=t-E_1-\ldots-E_n,\qquad\forall 0\leq n\leq k
$$
and
$$
\ell_{0(k+1)}=0\vee(t-E_1-\ldots-E_{k+1}),
$$
where $(E_n)_{n\geq 1}$ is a sequence of i.i.d.\ exponential r.v.\ of parameter $\theta$. In
addition, on the event $\{\ell_{0k}\geq s_1,\ \ell_{0k}-\ell_{0(k+1)}<s_1\}$, one has
$\ell_{0k}-\ell_{0(k+1)}=E_{k+1}$. Hence,
\begin{multline*}
  \PP_t (A_{0j}, A_{0k}, k\leq M_0)
  \leq\frac{b}{\alpha}\,\EE_t\Bigg[\mathbbm{1}_{L_{0j}\in[s_1,s_2]}\,\PP_t\left(R_t^{0j}\geq
    x\mid \ell_{0j}=L_{0j},\ell_{0(j+1)}=L_{0(j+1)}\right)\times \\
  \times \left.  \mathbbm{1}_{E_{j+2}+\ldots+E_k\leq s_2-s_1}\left(1-\mathbbm{1}_{E_{k+1}<s_1}
      \frac{W_\theta(s_1-E_{k+1})}{W_\theta(s_2)}\right)
    \left(1-\frac{1}{W_\theta(s_2)}\right)^{\lceil x\rceil-1}\right],
\end{multline*}
where
$$
L_{0j}=0\vee(t-E_1-\ldots-E_j)\qquad\text{and}\qquad L_{0(j+1)}=0\vee(t-E_1-\ldots-E_{j+1}).
$$
Since $(L_{0j},L_{0(j+1)})$,\ $(E_{j+2},\ldots,E_k)$ and $E_{k+1}$ are independent, we finally
obtain
\begin{align}
  \sum_{0\leq j<k}\PP_t (A_{0j}, A_{0k}, k\leq M_0)& \leq\frac{b}{\alpha}\sum_{0\leq
    j<k}\PP_t(A_{0j},j\leq M_0)\PP(E_{j+2}+\ldots+E_k\leq s_2-s_1)\times \notag \\ & \qquad
  \times\left(1-\EE\left(\mathbbm{1}_{E_{k+1}<s_1}\frac{W_\theta(s_1-E_{k+1})}{W_\theta(s_2)}\right)\right)
  \left(1-\frac{1}{W_\theta(s_2)}\right)^{\lceil x\rceil-1} \notag \\ &
  =\frac{b}{\alpha}\,(\EE_t K_t(x,s_1,s_2))\sum_{i\geq
    0}\PP(E_1+\ldots+E_i\leq s_2-s_1)\times \notag \\ & \qquad \times\left(1-\int_0^{s_1}\theta
    e^{-\theta z}\frac{W_\theta(s_1-z)}{W_\theta(s_2)}dz\right)
  \left(1-\frac{1}{W_\theta(s_2)}\right)^{\lceil x\rceil-1}. \label{eq:youpla-2}
\end{align}
Now, we have
\begin{equation*}
  \sum_{i\geq 0}\PP(E_1+\ldots+E_i\leq s_2-s_1)=1+\EE(P)=1+\theta(s_2-s_1),
\end{equation*}
where $P$ is a Poisson r.v.\ of parameter $\theta(s_2-s_1)$. Combining this equation
with~(\ref{eq:youpla-1}) and~(\ref{eq:youpla-2}) ends the proof
of~(\ref{eq:K(K-1)}).\hfill$\Box$
\bigskip

\paragraph{Acknowledgments.} This work was funded by project MANEGE `Mod\`eles
Al\'eatoires en \'Ecologie, G\'en\'etique et \'Evolution'
09-BLAN-0215 of ANR (French national research agency).

\end{document}